\documentclass[11pt]{article}
\usepackage{amsfonts}
\usepackage{mathrsfs}
\usepackage{amssymb,amsmath} 
 \allowdisplaybreaks

\catcode`!=11
\let\!int\int \def\int{\displaystyle\!int}
\let\!lim\lim \def\lim{\displaystyle\!lim}
\let\!sum\sum \def\sum{\displaystyle\!sum}
\let\!sup\sup \def\sup{\displaystyle\!sup}
\let\!inf\inf \def\inf{\displaystyle\!inf}
\let\!cap\cap \def\cap{\displaystyle\!cap}
\let\!max\max \def\max{\displaystyle\!max}
\let\!min\min \def\min{\displaystyle\!min}
\let\!frac\frac \def\frac{\displaystyle\!frac}
\catcode`!=12

\let\oldsection\section
\renewcommand\section{\setcounter{equation}{0}\oldsection}

\allowdisplaybreaks
\def\pf{\it{Proof.}\rm\quad}

\newcommand\uu{{\bf u}}
\newcommand\BB{{\bf B}}
\newcommand\xx{{\bf x}}
\newcommand\divg{{\text{div}}}
\newtheorem{defn}{Definition}[section]
\newtheorem{thm}{Theorem}[section]

\newtheorem{lem}{Lemma}[section]
\newtheorem{re}{Remark}[section]
\begin{document}
%
\title{\bf On the vanishing resistivity
limit and the magnetic boundary-layers for one-dimensional
compressible magnetohydrodynamics\thanks{Supported by the National
Basic Research Program under the grant 2011CB309705, NSFC (Grant
Nos. 11271306, 11229101, 11371065) and the Beijing Center for
Mathematics and Information Interdisciplinary Sciences.}}
\author{
{\bf Song Jiang$^1$\thanks{Email: jiang@iapcm.ac.cn},  \qquad Jianwen Zhang$^2$\thanks{Email: jwzhang@xmu.edu.cn (Corresponding author)}}
 \\[1mm]
\small$^1$ Institute of Applied Physics and Computational Mathematics,\\
\small P.O. Box 8009, Beijing 100088, China\\[1mm]
\small$^2$ School of Mathematical Sciences, Xiamen University,
Xiamen 361005, China}
\date{}
\maketitle

\noindent{\bf Abstract.} We consider an initial-boundary value
problem for the one-dimensional equations of compressible isentropic
viscous and non-resistive magnetohydrodynamic flows. The global
well-posedness of strong solutions with general large data is
established. Moreover, the vanishing resistivity limit is
justified and the thickness of magnetic boundary layers is analyzed.
The proofs of these results are based on a full use of the so-called
``effective viscous flux'', the material derivative and the
structure of the equations.

\vskip 2mm


\vskip 2mm


\section{Introduction}
Magnetohydrodynamics (MHD) concerns the motion of a conducting fluids in an
electromagnetic field and has a very wide range of applications in astrophysics,
plasma, and so on.  Because the dynamic
motion of the fluids and the magnetic field interacts strongly on each other,
both the hydrodynamic and electrodynamic effects have to be considered. The three-dimensional
equations for  compressible isentropic
magnetohydrodynamic flows, derived from fluid mechanics with
appropriate modifications to account for electrical forces, read as follows (see \cite{Ca,KL}):
\begin{equation}\label{1.1}
\begin{cases}
\rho_t+\divg(\rho \uu)=0,\\
(\rho \uu)_t+\divg(\rho \uu\otimes\uu)+\nabla P=\mu\Delta \uu+(\mu+\lambda')\nabla\divg\uu+(\nabla \times
\BB)\times
\BB,\\
\BB_t-\nabla\times(\uu\times\BB)=-\nu\nabla\times(\nabla\times\BB),\quad\divg\BB=0
\end{cases}
\end{equation}
with $\xx\in\Omega\subset\mathbb{R}^3$ and $t\geq0$. Here, the
unknown functions $\rho, \uu\in\mathbb{R}^3, P$ and
$\BB\in\mathbb{R}^3$ are the density of fluid, velocity, pressure,
and magnetic field, respectively.  The viscosity coefficients $\mu$
and $\lambda'$ satisfy the physical conditions $ \mu>0,
3\lambda'+2\mu\geq 0$. The constant $\nu>0$ is the resistivity
coefficient  which is inversely proportional to the electrical
conductivity constant (magnetic Reynolds number) and acts as the
diffusivity coefficient of the magnetic fields.  The pressure $P(\rho)$
is generally determined through the equation of state (the so-called
$\gamma$-law):
\begin{equation}
P(\rho)\triangleq A\rho^\gamma\quad{\rm with}\quad
A>0,\gamma>1.\label{1.2}
\end{equation}

Equations (\ref{1.1}), (\ref{1.2})  describe the macroscopic
behavior of electrically conducting compressible (isentropic) fluids
in a magnetic field. From Eq. (\ref{1.1})$_3$ it is clear that the
time rate of change of the magnetic field (i.e., $\BB_t$) is dominated
by both the advection term $\nabla\times(\uu\times\BB)$ and the
diffusion term $\nu\nabla\times(\nabla\times\BB)$. However, in many
cosmical and geophysical problems where the conducting fluid is of
extremely high conductivity (ideal conductors), the resistivity
$\nu$ is inversely proportional to the electrical conductivity
$\sigma$, and therefore, it is more rational to assume that there is no
magnetic diffusion (i.e. $\nu=0$) (see, for example,
\cite{Ch1961,Fr1982}). So, instead of (\ref{1.1})$_3$, the induction
equation for magnetic field in such cases has the form:
$$
\BB_t-\nabla\times(\uu\times\BB)=0,
$$
which particularly  implies  that in a highly conducting fluid the
magnetic field lines move along exactly with the fluid, rather than
simply diffusing out. This type of behavior is physically expressed
as that {\it the magnetic field lines are frozen into the fluid}. In
effect, the fluid can flow freely along the magnetic field lines,
but any motion of the conducting fluid, perpendicular to the field
lines, carries them with the fluid. The ``frozen-in" nature of
magnetic fields plays very important roles and has a very wide range
of applications in both astrophysics and nuclear fusion theory,
where the magnetic Reynolds number $R_m\thicksim 1/\nu$ is usually
very high. A typical illustration of the ``frozen-in" behavior is
the phenomenon of sunspots. For more details of its physical
background and applications, we refer to
\cite{Bi,BS2003,Ca,Ch1961,Fr1982,GB,KL,LLP}.

Formally, when $\nu=0$, system (\ref{1.1}) turns into
\begin{equation}\label{1.3}
\begin{cases}
\rho_t+\divg(\rho \uu)=0,\\
(\rho \uu)_t+\divg(\rho \uu\otimes\uu)+\nabla P=\mu\Delta \uu+(\mu+\lambda')\nabla\divg\uu+(\nabla \times
\BB)\times
\BB,\\
\BB_t-\nabla\times(\uu\times\BB)=0,\quad\divg\BB=0,
\end{cases}
\end{equation}
where the pressure $P(\rho)$ satisfies the $\gamma$-law
(\ref{1.2}). This is often called the compressible isentropic
viscous and non-resistive MHD equations.

Because of the complete interaction between the dynamic motion and
the magnetic field, the strong nonlinear terms and the lack of
dissipation mechanism of the magnetic field, many physically important
and mathematically fundamental problems of system (\ref{1.3}) are
still open. For example, to the author's knowledge, there is no any
rigorously mathematical literature on the global well-posedness of
the initial (boundary) value problem of (\ref{1.3}), even that the
initial data are close to a non-vacuum equilibrium, though the same
problem has been successfully solved for the compressible
Navier-Stokes equations (i.e., $\BB=0$) by Matsumura-Nishida
\cite{MN1980}.

Due to the complicated structure of multi-dimensional equations,
instead of (\ref{1.1}) and (\ref{1.3}), in this paper we shall
consider the simplest one-dimensional equations (cf. \cite{Gu1964}):
\begin{equation}
\begin{cases}
\rho_t+(\rho u)_x=0,\\
(\rho u)_t+\left(\rho
u^2+P(\rho)+\frac{1}{2}b^2\right)_x=\lambda
u_{xx},\\
b_t+(ub)_x=\nu b_{xx},
\end{cases}\label{1.4}
\end{equation}
and
\begin{equation}
\begin{cases}
\rho_t+(\rho u)_x=0,\\
(\rho u)_t+\left(\rho
u^2+P(\rho)+\frac{1}{2}b^2\right)_x=\lambda
u_{xx},\\
b_t+(ub)_x=0,
\end{cases}\label{1.5}
\end{equation}
where the pressure $P(\rho)$ obeys the $\gamma$-law (\ref{1.2}) and
$\lambda=2\mu+\lambda'$.

Indeed, assume that the MHD flow is moving only in the longitudinal
direction $x$ and uniform in the transverse directions  $(y,z)$.
Then it is easy to derive (\ref{1.4}) and (\ref{1.5}) from
(\ref{1.1}) and (\ref{1.3}), respectively, based on the specific
choice of dependent variables:
$$
\rho=\rho(x,t),\quad \uu=(u(x,t),0,0),\quad\BB =(0,0,b(x,t)).
$$
We mention here that the one-dimensional system for compressible
heat-conductive viscous and resistive MHD flows in a form similar to
that in (\ref{1.4}) was studied by Kazhikhov-Smagulov \cite{KS},
where the global well-posedness of solutions was announced.

The main purpose of this paper is to show the global well-posedness
of strong solutions to an initial-boundary value problem of the system
(\ref{1.5}), to justify the vanishing resistivity limit (i.e.,
$\nu\to0$) from (\ref{1.4}) to (\ref{1.5}) rigorously, and to study
the boundary layer effects as $\nu\to0$. To do so, without loss of
generality, we  consider an initial-boundary value problem of
(\ref{1.5}) on a bounded spatial-domain $\Omega\triangleq(0,1)$ with
the following initial and boundary data:
\begin{equation}
\begin{cases}
(\rho,u,b)(x,0)=(\rho_0,u_0,b_0)(x),\quad x\in[0,1],\\
u(0,t)=u(1,t)=0,\quad t\geq0.
\end{cases}\label{1.6}
\end{equation}

Our first result concerns the global well-posedness of strong
solutions to the initial-boundary value problem (\ref{1.5}),
(\ref{1.6}).
\begin{thm}\label{thm1.1}Assume that the initial data
$(\rho_0,u_0,b_0)$ given in (\ref{1.6}) satisfies
\begin{equation}
\inf_{0\leq x\leq1}\rho_0(x)>0,\quad (\rho_0,b_0)\in H^1,\quad
u_0\in H_0^1\cap H^2.
\end{equation}
Then for any $0<T<\infty$, there exists a global unique strong
solution $(\rho,u,b)$ to the initial-boundary value problem
(\ref{1.5}), (\ref{1.6}) on $[0,1]\times[0,T)$, satisfying
\begin{equation}
0<C^{-1}\leq \rho(x,t)\leq C<\infty,\quad \forall\
(x,t)\in[0,1]\times[0,T)\label{1.8}
\end{equation}
for some positive constant $C$, and
\begin{equation}\label{1.9}
\begin{cases}
(\rho,b)\in L^\infty(0,T;H^1),\quad (\rho_t,b_t)\in
L^\infty(0,T;L^2),\\
 u\in L^\infty(0,T;H_0^1\cap
H^2),\quad  u_t\in L^\infty(0,T;L^2)\cap L^2(0,T;H^1).
\end{cases}
\end{equation}
\end{thm}

Theorem \ref{thm1.1} is the first result on the global
well-posedness theory of the non-resistive MHD equations with large
data. It particularly implies that the solutions of the viscous and
non-resistive MHD problem (\ref{1.5}), (\ref{1.6}) will not develop
vacuum and mass concentration in finite time provided the
initial data are bounded and smooth, and do not contain vacuum.

It is worth pointing out that the equations (\ref{1.5}) look similar to
a compressible model for gas and liquid two-phase fluids (see, for example, \cite{YZ2011,Ev2011}).
However, to prove the global well-posedness of the two-phase model,
it is technically assumed in \cite{YZ2011,Ev2011} that the
proportion between the mass of gas and liquid must be bounded, in
analogy to the assumptions  $\rho,b\geq0$ and $0\leq b/\rho<\infty $
for (\ref{1.5}). Of course, this is somewhat reasonable for the
two-phase model, but not physical and realistic in
magnetohydrodynamics.

The second purpose of this paper is to justify the
vanishing resistivity limit from (\ref{1.4}) to (\ref{1.5})
rigorously, as the resistivity coefficient $\nu\to0$. More
precisely,
\begin{thm}\label{thm1.2}(i)
Suppose that system (\ref{1.4}) is equipped with initial and
boundary data:
\begin{equation}
\begin{cases}
(\rho,u,b)(x,0)=(\rho_0,u_0,b_0)(x),\quad x\in[0,1],\\
u(0,t)=u(1,t)=0,\quad b(0,t)= b_1(t),\quad b(1,t)=b_2(t),\quad
t\geq0.
\end{cases}\label{1.10}
\end{equation}
For any $0<T<\infty$, assume that
\begin{equation}\label{1.11}
\inf_{0\leq x\leq1}\rho_0(x)>0,\quad (\rho_0,b_0)\in H^1,\quad
(b_1,b_2)\in C^1([0,T)),\quad u_0\in H_0^1.
\end{equation}
Then the initial-boundary value problem (\ref{1.4}), (\ref{1.10})
has a global unique strong solution $(\rho,u,b)$ on
$[0,1]\times[0,T)$, satisfying
\begin{equation}
0<C^{-1}\leq \rho(x,t)\leq C<\infty,\quad \forall\
(x,t)\in[0,1]\times[0,T)\label{1.12}
\end{equation}
 and
\begin{eqnarray}
&&\sup_{0\leq t< T}\left(\|u_x \|_{L^2}^2+\|b
\|_{L^\infty}^2+\nu^{1/2}\|\rho_x\|_{L^2}^2+\nu^{1/2}\|b_x\|_{L^2}^2\right)(t)\nonumber\\
&&\qquad +\int_0^T\left(\|\rho^{1/2}\dot
u\|_{L^2}^2+\nu^{3/2}\|b_{xx}\|_{L^2}^2+\nu^{1/2}\|b
b_{x}\|_{L^2}^2\right)dt\leq C,\label{1.13}
\end{eqnarray}
where $C$ is a positive constant independent of $\nu$.

\vskip 2mm

(ii) Assume that $(\rho^\nu,u^\nu,b^\nu)$ and $(\rho,u,b)$, defined
on $[0,1]\times [0,T)$, are the solutions of the problems
(\ref{1.4}), (\ref{1.10}) and (\ref{1.5}), (\ref{1.6}),
respectively. Then,
\begin{equation*}
\begin{cases}
(\rho^\nu,u^\nu,b^\nu)\to (\rho,u,b)\quad {strongly\ \ in}\quad L^\infty(0,T;L^2),
\\
\nu b_x^\nu\to0,\quad u^\nu_x\to u_x \quad {strongly\ \ in}\quad L^2(0,T;L^2),
\end{cases}
\end{equation*}
and moreover, there exists a positive constant $C$, independent
of $\nu$, such that for $\nu\in(0,1)$,
\begin{equation}\label{1.14}
\sup_{0\leq t<T}\left(\|\rho^\nu-\rho\|_{L^2}^2+
\|b^\nu-b\|_{L^2}^2+\| u^\nu-u\|_{L^2}^2\right)(t) +\int_0^T
\|(u^\nu-u)_x\|_{L^2}^2 dt\leq C\nu^{1/2}.
\end{equation}
\end{thm}

The third and main result of this paper is to study the effects of
magnetic boundary layer as the resistivity coefficient $\nu\to0$. In
fact, it is obvious that when the resistivity coefficient goes to
zero, the parabolic equation (\ref{1.4})$_3$ turns into the
hyperbolic equation (\ref{1.5})$_3$, and the boundaries become
characteristic due to the non-slip boundary conditions
$u|_{x=0,1}=0$. Thus, by the classical theory in \cite{OR1973}, one
has to drop the boundary conditions of the magnetic field in
(\ref{1.10}) (cp. (\ref{1.6})), and consequently, because of the
disparity of boundary conditions, we cannot expect that as
$\nu\to0$, the solution of the problem (\ref{1.4}), (\ref{1.10})
will tend to the one of the problem (\ref{1.5}), (\ref{1.6})
uniformly up to the boundaries $x=0,1$. In other words, a (magnetic)
boundary layer appears near the boundary.

Similar to the relations among the Euler, Navier-Stokes and Prandtl equations (see, for example,
\cite{LL1987,OS1999,Sc1987}), it is expected that as $\nu\to0$, the solution of the problem (\ref{1.4}), (\ref{1.10}) converges uniformly
to the solution of the problem ({\ref{1.5}), (\ref{1.6}) away from the boundaries, while there is a sharp change of
gradient near the boundary. Inspired by this, we introduce the concept of magnetic boundary-layer thickness
(MBL-thickness) as follows.
\begin{defn}\label{def1.1}
A non-negative function $\delta(\nu)$ is called a MBL-thickness of the problem (\ref{1.4}), (\ref{1.10})
with vanishing resistivity limit, if $\delta(\nu)\downarrow0$ as $\nu\downarrow0$, and
\begin{equation}
\lim_{\nu\to0}\left\|\left(\rho^\nu-\rho,u^\nu-u,b^\nu-b\right)\right\|_{L^\infty(0,T;C(\overline\Omega_{\delta(\nu)}))}=0,\label{1.15}
\end{equation}
\begin{equation}
\liminf_{\nu\to0}\left\|\left(\rho^\nu-\rho,u^\nu-u,b^\nu-b\right)\right\|_{L^\infty(0,T;C(\overline\Omega))}>0
\label{1.16}
\end{equation}
with $\Omega_\delta\triangleq\{x\in\Omega\ |\ \delta<x<1-\delta\}$. Here, $(\rho^\nu,u^\nu,b^\nu)$
and $(\rho, u,b)$ are the solutions of the problems
(\ref{1.4}), (\ref{1.10}) and (\ref{1.5}), (\ref{1.6}), respectively.
\end{defn}

The concept of boundary-layer thickness (BL-thickness)  has been
introduced in \cite{FS1999,JZ2009} in a similar manner for the
one-dimensional cylindrical compressible Navier-Stokes equations
with vanishing shear viscosity limit. The BL-thickness for the
scalar conservation laws and the 2D Boussinesq equations with
vanishing diffusivity limit was also studied in \cite{FS2004} and
\cite{JZZ2011}, respectively. It is worth mentioning that Definition
\ref{def1.1} does not determine the MBL-thickness uniquely, since
any function $\delta_1(\nu)$, satisfying $\delta_1(\nu)\geq
\delta(\nu)$ and $\delta_1(\nu)\downarrow 0$ as $\nu\downarrow0$, is
also a MBL-thickness. Thus, there should exist a minimal
MBL-thickness $\delta_*(\nu)$ which may be considered as the true
MBL-thickness.

In this paper, we shall prove that a function
$\delta_n(\nu)=\nu^{1/2-1/n}$ with $n>2$ is a MBL-thickness in the
sense of Definition \ref{def1.1}. This is somewhat in agreement with
the famous Stokes-Blasius law in the laminar boundary layer theory
(see \cite{Sc1987}), since
$\liminf_{n\to\infty}\delta_n(\nu)=\nu^{1/2}$. In order to simplify
the analysis of MBL-thickness, similarly to \cite{FS1999,JZ2009}, we shall focus on the special case of
vanishing initial data. More precisely, we shall prove
\begin{thm}\label{thm1.3} (i) Assume that the viscous and non-resistive system (\ref{1.5}) is equipped with initial and boundary data:
\begin{equation}
(\rho,u,b)|_{t=0}=(\overline\rho,0,0)\quad {with}\quad \overline\rho\equiv{Const.}>0,\quad {and}\quad u|_{x=0,1}=0.\label{1.17}
\end{equation}
Then the problem (\ref{1.5}), (\ref{1.17}) has only a trivial solution $(\rho,u,b)=(\overline\rho,0,0)$.

\vskip 2mm

(ii) Assume that $(\rho,u,b)$ is the solution of the viscous and resistive system (\ref{1.4}) with the following initial and boundary data:
\begin{equation}
(\rho,u,b)|_{t=0}=(\overline\rho,0,0),\quad u|_{x=0,1}=0,\quad b(0,t)=b_1(t),\quad b(1,t)=b_2(t),\label{1.18}
\end{equation}
where $\overline\rho\equiv Const.>0$ and $b_1(t),b_2(t)$ are the same as in (\ref{1.17}) and (\ref{1.10}), respectively.
Then any function $\delta (\nu)\geq 0$, satisfying
\begin{equation}
\delta(\nu)\to 0\quad{and}\quad \frac{\delta(\nu)}{\nu^{1/2}}\to\infty,\quad{as}\quad \nu\to0,\label{1.19}
\end{equation}
is a MBL-thickness in the sense of Definition \ref{def1.1} such that
\begin{equation}
\lim_{\nu\to0}\left(\|\rho-\overline\rho\|_{L^\infty(0,T;C(\overline\Omega_{\delta(\nu)}))}^2+\|b\|_{L^\infty(0,T;C(\overline\Omega_{\delta(\nu)}))}^2\right)=0\label{1.20}
\end{equation}
and
\begin{equation}
\liminf_{\nu\to0}\left(\|\rho-\overline\rho\|_{L^\infty(0,T;C(\overline{\Omega}))}^2+\|b\|_{L^\infty(0,T;C(\overline\Omega))}^2\right)>0,\label{1.21}
\end{equation}
provided the boundary data $b_1(t),b_2(t)$ are not identically zero.
\end{thm}

In view of Definition \ref{def1.1}, we notice from (\ref{1.20}) and
(\ref{1.21})  that there is no boundary layer effect on the
velocity. This is mainly due to the smooth  mechanism of the
viscosity term $\lambda u_{xx}$ with $\lambda>0$, and it indeed
holds that (see (\ref{4.4}), (\ref{4.5}) below)
$$
\|u\|_{L^\infty(0,T;C(\overline\Omega))}\to0,\quad{\rm as}\quad \nu\to0.
$$

The proofs of Theorems \ref{thm1.1}--\ref{thm1.3} will be given respectively in Sections 1, 2 and 3, based on the
global (uniform-in-$\nu$) a-priori estimates of the solutions.
However, the lack of smooth mechanism of the magnetic field, the strong nonlinearities, and the
interaction between  dynamic motion and magnetic field
will cause some serious difficulties. To circumvent these difficulties,  motivated by the study of multi-dimensional
Navier-Stokes/MHD equations (see, for example,
\cite{FNP,Li,LXZ2013}), we shall make
a full use of the so-called ``effective viscous flux" $F$ and the
material derivative $\dot u$:
\begin{equation}
F\triangleq \lambda u_x-P(\rho)-\frac{b^2}{2}\quad {\rm and}\quad \dot u\triangleq u_t+uu_x.\label{1.22}
\end{equation}
It turns out that the ``effective viscous flux" $F$ possesses more
regularities than the velocity.

We begin the proofs with the elementary energy estimates and the upper boundedness of the density (see Lemmas \ref{lem2.1},
\ref{lem2.2}, \ref{lem3.2} and \ref{lem3.3}). Next, by using the
special mathematical structure of $F$, especially, using the
non-negativity of $b^2$ and the fact that $F_x=\rho\dot
u$ due to (\ref{1.4})$_2$ (or (\ref{1.5})$_2$), we can improve the
integrability of the magnetic field and obtain the desired bounds of
$\|u_x\|_{L^\infty(0,T;L^2)}$ and $\|\rho^{1/2}\dot
u\|_{L^2(0,T;L^2)}$ (see Lemmas \ref{lem2.3}--\ref{lem2.5} and
\ref{lem3.4}). With these estimates at hand, we then can show the
boundedness of the magnetic field and the lower boundedness of the density as well
(see Lemmas \ref{lem2.6} and \ref{lem3.5}). Thus, noting that
$F_x=\rho\dot u$,
  by (\ref{1.22}) we deduce
\begin{eqnarray}
\|u_x\|_{L^\infty}&\leq&C\left(\|F\|_{L^\infty}+\|P(\rho)\|_{L^\infty}+\|b^2\|_{L^\infty}\right)
\leq C\left(1+\|F\|_{L^\infty}\right)\nonumber\\
&\leq& C\left(1+\|F\|_{L^2}+\|F_x\|_{L^2}\right)\leq
C\left(1+\|\rho^{1/2}\dot u\|_{L^2}\right)\in L^2(0,T),\label{1.23}
\end{eqnarray}
which plays a very important role in the entire analysis of this paper.
Indeed, as an immediate result of (\ref{1.23})  and the blowup criterion in
\cite{XZ2011}, one can easily obtain the global existence of strong
solutions to (\ref{1.5}), (\ref{1.6}) (see Section 2).

The justification of vanishing resistivity limit and the study of
magnetic boundary layers are more difficult and need some more
delicate estimates. Indeed, due to the presence of boundary layer
effects, the global uniform-in-$\nu$ estimates for the solution
$(\rho,u,b)$ of the problem (\ref{1.4}), (\ref{1.10}) (also
(\ref{1.4}), (\ref{1.18})) are much fewer than those for the solutions of the
problem (\ref{1.5}), (\ref{1.6}). For example, it is very difficult to
obtain the global uniform (in $\nu$) $L^2$-bounds of the derivatives
of the density and magnetic field for the problems (\ref{1.4}), (\ref{1.10})
and (\ref{1.4}), (\ref{1.18}), compared with the uniform bounds stated in
Lemma \ref{lem2.7} for the problem (\ref{1.5}), (\ref{1.6}). Instead of (\ref{2.27}), noting that
$$
u_{xx}=\lambda^{-1}\left(F_x+P(\rho)_x+bb_x\right),
$$
using (\ref{1.23}), the estimates obtained and subtle boundary
analysis, we have (cf. Lemma \ref{lem3.6})
\begin{equation}
\nu^{1/2}\sup_{0\leq t<
T}\left(\|\rho_x(t)\|_{L^2}^2+\|b_x(t)\|_{L^2}^2\right)+\nu^{3/2}\int_0^T\|b_{xx}\|_{L^2}^2dt\leq
C.\label{1.24}
\end{equation}
This suffices to prove the vanishing resistivity limit and to obtain
the convergence rates given in Theorem \ref{thm1.2} (see Section 3).

As for the analysis of boundary layer effects, the
weighted $L^1$-method used  in \cite{FS1999,JZ2009} seems difficult
to apply here, due to the strong interactions of dynamic
motion and magnetic field. Instead of the $L^1$-method, we shall make use of the weighted $L^2$-method to analyze the thickness of magnetic boundary layers.
To do this, we first utilize the special
initial data in (\ref{1.18}), (\ref{1.14}) and (\ref{1.24}) to  improve the convergence estimate of the velocity (cf. Lemma \ref{lem4.1}):
\begin{equation}
\sup_{0\leq t< T}\|u_x(t)\|_{L^2}^2+\int_0^T\|\rho^{1/2}\dot
u\|_{L^2}^2dt\leq C\nu^{1/2}.\label{1.25}
\end{equation}
This particularly implies that  $\|u(t)\|_{C(\overline\Omega)}\to 0$ as
$\nu\to 0$ for $t\in[0, T)$, and hence, there is no boundary layer
effects on the velocity. Using (\ref{1.22}),  (\ref{1.24}),
(\ref{1.25}) and the special initial data in (\ref{1.18}) again,
we succeed in deriving the weighted (interior) $L^2$-estimate (cf. Lemma \ref{lem4.2}):
\begin{equation}
\sup_{0\leq t<T}\int_0^1 \xi(x)\left(| \rho_x|^2+|b_x|^2\right)(x,t)
dx\leq C\nu^{1/2}\quad{\rm with}\quad \xi(x)\triangleq x^2(1-x)^2, \label{1.26}
\end{equation}
which, together with (\ref{1.14}) and Sobolev's inequality, proves Theorem \ref{thm1.3} (see Section 4).

\section{Global well-posedness of (\ref{1.5}), (\ref{1.6})}
The local existence of regular solutions with smooth initial data
can be shown by the standard method based on the Banach theorem and
the contractivity of the operator defined by the linearized the problem on a small time interval.
The global existence of solutions will be obtained by the method of extending local
solutions with respect to time based on the global a-priori estimates.

For this purpose, let $(\rho,u,b)$ be a smooth solution of
(\ref{1.5}), (\ref{1.6}). We will establish the necessary global
a-priori estimates of $(\rho,u,b)$ defined on $[0,1]\times[0,T)$ for
any fixed $T>0$. For simplicity, we denote by $C$ and $C_i$
$(i=1,2,\ldots)$ generic positive constants which may depend on
$\lambda,A,\gamma$, the initial norms of $(\rho_0,u_0,b_0)$ and $T$,
and may change from line to line.

We begin with the conservations of mass and momentum.
\begin{lem}\label{lem2.1} Let $(\rho,u,b)$ be a smooth solution of
(\ref{1.5}), (\ref{1.6}) on $[0,1]\times[0,T)$. Then for any $0\leq
t< T$,
\begin{equation}
0<\int_0^1\rho(x,t)dx=\int_0^1\rho_0(x)dx<\infty,\label{2.1}
\end{equation}
and
\begin{eqnarray}
&&\int_0^1\left(\frac{1}{2}\rho
u^2+\frac{1}{2}b^2+\frac{A}{\gamma-1}\rho^\gamma\right)(x,t)dx+\lambda\int_0^t\|
u_x\|_{L^2}^2 ds\nonumber\\
&&\qquad =\int_0^1\left(\frac{1}{2}\rho_0
u_0^2+\frac{1}{2}b_0^2+\frac{A}{\gamma-1}\rho_0^\gamma\right)(x)dx.\label{2.2}
\end{eqnarray}
\end{lem}

The upper boundedness of the density can be deduced in a similar manner to that in \cite{FJN2007}.
\begin{lem}\label{lem2.2}Let $(\rho,u,b)$ be a smooth solution of
(\ref{1.5}), (\ref{1.6}) on $[0,1]\times[0,T)$. Then,
\begin{equation}
0\leq \rho(x,t) \leq C,\quad\forall\
(x,t)\in[0,1]\times[0,T).\label{2.3}
\end{equation}
\end{lem}
\pf The non-negativity of the density (i.e., $\rho\geq0$) readily
follows from the method of characteristics and the fact that $\rho_0>0$. In the next,
for completeness we sketch the proof of the upper bound below. Define
\begin{equation}
\psi(x,t)\triangleq\int_0^t\left(\lambda u_x-\rho u^2
-P(\rho)-\frac{b^2}{2}\right)(x,s)ds+\int_0^x(\rho_0u_0)(\xi)d\xi.\label{2.4}
\end{equation}

Clearly, it follows from (\ref{1.5})$_1$ and (\ref{1.5})$_2$ that
\begin{equation}\label{2.5}
\psi_x=\rho u,\quad\psi_t=\lambda u_x-\rho u^2-P(\rho)
-\frac{b^2}{2},\quad\psi|_{t=0}=\int_0^x(\rho_0u_0)(\xi)d\xi,
\end{equation}
and hence, by Lemma \ref{lem2.1} one has
$$
\left|\int_0^1 \psi(x,t)dx\right|\leq C\quad{\rm
and}\quad\|\psi_x\|_{L^\infty(0,T;L^1)}\leq C,
$$
which particularly yields
\begin{equation}
\|\psi\|_{L^\infty(0,T;L^\infty)}\leq
\left|\int_0^1\psi(x,t)dx\right|+\|\psi_x\|_{L^\infty (0,T;L^1)}
\leq C.  \label{2.6}
\end{equation}

Let $D_t\triangleq\partial_t+u\partial_x$ denote the material
derivative and set
\begin{equation}\label{2.7}
\Phi(x,t)\triangleq\exp\left\{\frac{\psi(x,t)}{\lambda}\right\}.
\end{equation}
Then, by straightforward calculations we have
$$
D_t(\rho \Phi)=\partial_t(\rho \Phi)+u\partial_x(\rho
\Phi)=-\frac{1}{\lambda}\left(P(\rho)+\frac{b^2}{2}\right)\rho
\Phi\leq0,
$$
and consequently,
$$
\|(\rho\Phi)(t)\|_{L^\infty}\leq \|(\rho\Phi)(0)\|_{L^\infty}\leq C,
$$
which, combined with the fact that $C^{-1}\leq \Phi(x,t)\leq C$ for
all $(x,t)\in[0,1]\times[0,T)$ due to (\ref{2.6}), immediately leads
to (\ref{2.3}).\hfill$\square$

\vskip 2mm

Due to lack of dissipation mechanism of the magnetic field and
the strong coupling of dynamic motion and magnetic field, the lower
boundedness of the density, which relies strongly on the $L^\infty$-norm of
magnetic field $b$, is more difficult to achieve, compared with the
upper bound. To circumvent the difficulties, motivated by the
mathematical theory of multi-dimensional Navier-Stokes/MHD equations
(see, for example, \cite{FNP,Li,LXZ2013}), we introduce the
so-called ``effective viscous flux" $F$, which possesses more
regularities than the velocity. Define
\begin{equation}\label{2.8}
F(x,t)\triangleq\left(\lambda u_x-P(\rho)-\frac{b^2}{2}
\right)(x,t)\quad {\rm and}\quad \dot u(x,t)\triangleq
(u_t+uu_x)(x,t),
\end{equation}
where $F$ and `` $\dot{}$ '' are the ``effective viscous flux" and
the material derivative, respectively. Thanks to (\ref{1.5})$_1$ and (\ref{1.5})$_2$, it is easy to see that
\begin{equation}
F_x=\rho (u_t+ u u_x)=\rho\dot u.\label{2.9}
\end{equation}
The quantities of ``effective viscous flux" $F$ and material
derivative $\dot u$ will play an important role in the entire
analysis, particularly, in controlling the first-order derivatives of
the solutions and studying the vanishing resistivity limit and the
boundary layer effects.

In order to estimate the first-order derivative of the velocity, we
first need the following preliminary lemma.
\begin{lem}\label{lem2.3} Let $F$  and `` $\dot{}$ '' be the same as in (\ref{2.8}). Then for any
$0\leq t< T$,
\begin{equation}
\frac{\lambda}{2}\frac{d}{dt}\|u_x\|_{L^2}^2+\|\rho^{1/2} \dot
u\|_{L^2}^2  \leq \frac{d}{dt}\int_0^1 \left(P(\rho)+\frac{b^2}{2}
\right) u_x dx+
C_1\left(\|F\|_{L^3}^3+\|b\|_{L^6}^6\right)+C.\label{2.10}
\end{equation}
\end{lem}
\pf Multiplying (\ref{1.5})$_2$ by $\dot u$ in $L^2$ and integrating
by parts, we deduce
\begin{eqnarray}
\frac{\lambda}{2}\frac{d}{dt}\int_0^1u_x^2dx+\int_0^1\rho \dot u^2 dx
&=&-\int_0^1 P(\rho)_x(u_t+uu_x) dx\nonumber\\
&&-\frac{1}{2}\int_0^1
(b^2)_x(u_t+uu_x)dx -\frac{\lambda}{2}\int_0^1 u_x^3 dx.\label{2.11}
\end{eqnarray}

Due to (\ref{1.5})$_1$, it holds that
\begin{equation}
P(\rho)_t+u P(\rho)_x+\gamma P(\rho) u_x=0\quad{\rm with}\quad
P(\rho)=A\rho^\gamma,\label{2.12}
\end{equation}
and hence, the first term on the right-hand side of (\ref{2.11}) can be estimated as follows.
\begin{eqnarray}
&&-\int_0^1 P(\rho)_x(u_t+uu_x) dx=\int_0^1
\left(P(\rho) u_{xt}-P(\rho)_x u u_x\right) dx\nonumber\\
&&\quad=\frac{d}{dt}\int_0^1 P(\rho) u_x dx-\int_0^1\left( P(\rho)_t
u_x+P(\rho)_x u u_x\right)dx\nonumber\\
&&\quad=\frac{d}{dt}\int_0^1 P(\rho) u_x dx+\gamma\int_0^1  P(\rho)
u_x^2 dx.\label{2.13}
\end{eqnarray}

Similarly, due to (\ref{1.5})$_3$, one has
\begin{equation}
(b^2)_t+u(b^2)_x+2b^2u_x=0,\label{2.14}
\end{equation}
so that,
\begin{eqnarray}
&&-\frac{1}{2}\int_0^1 (b^2)_x(u_t+uu_x)dx=
\frac{1}{2}\int_0^1\left( b^2
u_{xt}-(b^2)_x u u_x\right)dx\nonumber\\
&&\quad= \frac{1}{2}\frac{d}{dt}\int_0^1  b^2 u_{x}
dx-\frac{1}{2}\int_0^1\left((b^2)_t u_x+(b^2)_x u u_x\right)dx\nonumber\\
&&\quad= \frac{1}{2}\frac{d}{dt}\int_0^1  b^2 u_{x} dx+ \int_0^1
b^2u_x^2 dx.\label{2.15}
\end{eqnarray}

By virtue of (\ref{2.3}), we infer from (\ref{2.8}) that
\begin{equation}
|u_x|^3\leq C(1+|F|^3+b^6),\quad\forall\
(x,t)\in[0,1]\times[0,T).\label{2.16}
\end{equation}
Therefore, substituting (\ref{2.13}), (\ref{2.15}) into
(\ref{2.11}), and using (\ref{2.3}), (\ref{2.16}) and Cauchy-Schwarz's
inequality, we immediately obtain (\ref{2.10}).\hfill$\square$

\vskip 2mm

The next lemma is concerned with the higher integrability of the magnetic field.
\begin{lem}\label{lem2.4}Let $F$  and `` $\dot{}$ '' be the same ones defined in (\ref{2.8}). Then for any
$0\leq t< T$,
\begin{eqnarray}
\frac{d}{dt}\|b\|_{L^4}^4+ \|b\|_{L^6}^6\leq
C+C\|F\|_{L^3}^3.\label{2.17}
\end{eqnarray}
\end{lem}
\pf Indeed, multiplying (\ref{1.5})$_3$ by $4b^3$ and integrating by
parts, we get from (\ref{2.8}) that
\begin{eqnarray*}
0&=&\frac{d}{dt}\int_0^1 b^4 dx+3\int_0^1 u_x b^4 dx\nonumber\\
&=&\frac{d}{dt}\int_0^1 b^4
dx+\frac{3}{\lambda}\int_0^1\left(F+P(\rho)+\frac{ b^2}{2}\right)
b^4dx\nonumber\\
&=&\frac{d}{dt}\int_0^1 b^4 dx+ \frac{3}{2\lambda}\int_0^1 b^6 dx+\frac{3}{\lambda}\int_0^1\left(F+P(\rho) \right)
b^4dx,
\end{eqnarray*}
which, combined with (\ref{2.3}) and Cauchy-Schwarz's inequality, leads to (\ref{2.17}).\hfill$\square$

\vskip 2mm

Combining (\ref{2.10}) with (\ref{2.17}), we easily obtain
\begin{lem}\label{lem2.5} Let $(\rho,u,b)$ be a smooth solution of
(\ref{1.5}), (\ref{1.6}) on $[0,1]\times[0,T)$. Then,
\begin{equation}
\sup_{0\leq t<
T}\left(\|u_x(t)\|_{L^2}^2+\|b(t)\|_{L^4}^4\right)+\int_0^T\left(\|\rho^{1/2}\dot
u\|_{L^2}^2+\|b\|_{L^6}^6\right)dt\leq C,\label{2.18}
\end{equation}
and moreover,
\begin{equation}
\sup_{0\leq t<T}\|F(t)\|_{L^2}^2+\int_0^T\|F_x\|_{L^2}^2dt\leq
C,\label{2.19}
\end{equation}
where $F$ and `` $\dot{}$ " are defined in (\ref{2.8}).
\end{lem}
\pf Multiplying (\ref{2.17}) by a (large) number $K\geq \max\{(2\lambda)^{-1},2C_1\}$ and adding the
resulting inequality to (\ref{2.10}), we deduce
\begin{eqnarray*}
&&\frac{d}{dt}\int_0^1\left(\frac{\lambda}{4}u_x^2+\left(K-(4\lambda)^{-1}\right)b^4\right)dx+\int_0^1\left(\rho\dot
u^2+(K-C_1)b^6\right)dx\nonumber\\
&&\quad\leq \frac{d}{dt}\int_0^1P(\rho)u_x dx-\frac{1}{4}\frac{d}{dt}\int_0^1\left(\sqrt\lambda u_x-\frac{b^2}{\sqrt\lambda}\right)^2dx+C\int_0^1|F|^3dx+C,
\end{eqnarray*}
which, integrated over $(0,t)$, gives
\begin{equation}
 \left(\|u_x(t)\|_{L^2}^2+\|b(t)\|_{L^4}^4\right)+\int_0^t\left(\|\rho^{1/2}\dot
u\|_{L^2}^2+\|b\|_{L^6}^6\right)ds\leq
C+C\int_0^t\|F\|_{L^3}^3ds,\label{2.20}
\end{equation}
since integration of the second term on the right-hand side is
non-negative and it follows from  (\ref{2.3}) and Cauchy-Schwarz's inequality that
$$
\left|\int_0^1P(\rho)u_x dx\right|\leq
\frac{\lambda}{8}\|u_x\|_{L^2}^2+C\|P(\rho)\|_{L^2}^2\leq
\frac{\lambda}{8}\|u_x\|_{L^2}^2+C.
$$

Thanks to Lemmas \ref{lem2.1} and \ref{lem2.2}, we infer from
(\ref{2.8}) and (\ref{2.9}) that
\begin{equation}
\begin{cases}
\|F\|_{L^1}\leq
C\left(\|u_x\|_{L^1}+\|\rho\|_{L^\gamma}^\gamma+\|b\|_{L^2}^2\right)\leq
C\left(1+\|u_x\|_{L^2}\right),
\\[2mm]
\label{2.21} \|F\|_{L^2}^2\leq
C\left(1+\|u_x\|_{L^2}^2+\|b\|_{L^4}^4\right),\quad
\|F_x\|_{L^2}^2\leq C_2 \|\rho^{1/2} \dot u\|_{L^2}^2.
\end{cases}
\end{equation}
This, together with Sobolev's and Cauchy-Schwarz's inequalities, yields
($0<\varepsilon<1$)
\begin{eqnarray}
\|F\|_{L^3}^3&\leq& \|F\|_{L^1}\|F\|_{L^\infty}^2\leq
C(1+\|u_x\|_{L^2})\left(\|F\|_{L^2}^2+\|F\|_{L^2}\|F_x\|_{L^2}\right)
\nonumber\\
&\leq& C_2^{-1}
\varepsilon\|F_x\|_{L^2}^2+C\varepsilon^{-1}(1+\|u_x\|_{L^2}^2)
\|F\|_{L^2}^2 \nonumber\\
&\leq& \varepsilon\|\rho^{1/2} \dot
u\|_{L^2}^2+C\varepsilon^{-1}(1+\|u_x\|_{L^2}^2)\left(1
+\|u_x\|_{L^2}^2+\|b\|_{L^4}^4\right).\label{2.22}
\end{eqnarray}

Therefore, substituting (\ref{2.22}) into (\ref{2.20}), and choosing
$\varepsilon>0$ sufficiently small, by Gronwall's inequality we
obtain (\ref{2.18}), since (\ref{2.2}) implies that
$\|u_x(t)\|_{L^2}^2\in L^1(0,T)$. As a result,  (\ref{2.19}) follows
from (\ref{2.21}) and (\ref{2.18}).\hfill$\square$

\vskip 2mm

Next, by making a full use of the mathematical structure of  ``effective
viscous flux" $F$ again, we can prove the global boundedness of the
magnetic field, and consequently, the lower boundedness of the density can be
derived in a manner similar to that used in the derivation of the upper boundedness in Lemma \ref{lem2.2}.
\begin{lem}\label{lem2.6}Let $(\rho,u,b)$ be a smooth solution of (\ref{1.5}), (\ref{1.6}) on
$[0,1]\times[0,T)$. Then,
\begin{equation}
\sup_{0\leq t<
T}\left(\|b(t)\|_{L^\infty}+\|\rho^{-1}(t)\|_{L^\infty}\right)\leq
C.\label{2.23}
\end{equation}
\end{lem}
\pf Multiplying (\ref{1.5})$_3$ by $2n b^{2n-1}$ with $1\leq n\in \mathbb{N}$, integrating
by parts over $(0,1)$, and recalling the non-negativity of $P(\rho)$ and $b^2$, we obtain
\begin{eqnarray}
\frac{d}{dt}\int_0^1b^{2n}(x,t)dx&=&-(2n-1)\int_0^1 u_x
b^{2n}dx\nonumber\\
&=&-\frac{2n-1}{\lambda}\int_0^1
\left(F+P(\rho)+\frac{b^2}{2}\right)
b^{2n}dx\nonumber\\
&\leq&-\frac{2n-1}{\lambda}\int_0^1 F b^{2n}dx\leq C
n\|F\|_{L^\infty}\|b\|_{L^{2n}}^{2n}.\label{2.24}
\end{eqnarray}

It follows from (\ref{2.3}), (\ref{2.9}), (\ref{2.18}), (\ref{2.19}) and Sobolev's inequality that
\begin{equation}
\|F\|_{L^\infty}^2\leq
C\left(\|F\|_{L^2}^2+\|F\|_{L^2}\|F_x\|_{L^2}\right)\leq
C\left(1+\|\rho^{1/2}\dot u\|_{L^2}\right)\in L^2(0,T),\label{2.25}
\end{equation}
and hence, by direct calculations we deduce from (\ref{2.24}) and
(\ref{2.25}) that for $\forall$ $t\in[0,T]$,
\begin{equation*}
\|b(t)\|_{L^{2n}}^{2n}\leq C^{2n}
 \exp\left\{Cn\int_0^T\|F\|_{L^\infty}dt\right\}\leq
 C^{2n}\exp\left\{Cn\right\},
\end{equation*}
where $C$ is a positive constant independent of $n$. Thus, if  we
raise to the power $1/(2n)$ to both sides and let $n\to\infty$, then
we get
\begin{equation}
\|b(t)\|_{L^\infty}\leq C,\quad\forall\ t\in[0,T).\label{2.26}
\end{equation}

The lower boundedness of the density can be shown in the same way as in
Lemma \ref{lem2.2}. Let $\psi$ and $\Phi$ be the same functions
defined in (\ref{2.4}) and (\ref{2.7}), respectively. Analogously to
the proof of Lemma \ref{lem2.2}, we find
$$
D_t\left(\frac{1}{\rho\Phi}\right)=\frac{1}{\lambda}\left(P(\rho)+
\frac{b^2}{2}\right)\frac{1}{\rho \Phi}\quad{\rm with}\quad\emph{}
D_t\triangleq\partial_t+u\partial_x,
$$
so that, using (\ref{2.3}), (\ref{2.6}) and (\ref{2.26}), we have
\begin{equation*}
\left\|(\rho \Phi)^{-1}(t)\right\|_{L^\infty}\leq \left\|(\rho
\Phi)^{-1}(0)\right\|_{L^\infty}\exp\left\{\frac{1}{\lambda}\int_0^T\left(\|P(\rho)\|_{L^\infty}+\frac{1}{2}\|b\|_{L^\infty}^2\right)dt\right\}\leq
C,
\end{equation*}
which, combined with (\ref{2.6}) again, yields a desired lower bound
of the density.\hfill$\square$

\vskip 2mm

The next lemma concerns the first-order derivatives of the density and
magnetic field.
\begin{lem}\label{lem2.7}Let $(\rho,u,b)$ be a smooth solution of (\ref{1.5}), (\ref{1.6}) on
$[0,1]\times[0,T)$. Then,
\begin{equation}
\sup_{0\leq t< T}\left(\|\rho_x \|_{L^2}+\|b_x \|_{L^2}+\|\rho_t
\|_{L^2}+\|b_t \|_{L^2}+\|u \|_{L^\infty}\right)(t)\leq
C.\label{2.27}
\end{equation}
\end{lem}
\pf First, in view of (\ref{2.3}), (\ref{2.8}), (\ref{2.19}),
(\ref{2.23}) and Sobolev's inequality, we have
\begin{eqnarray}
\|u_x\|_{L^\infty} &\leq&
C\left(\|F\|_{L^\infty}+\|P(\rho)\|_{L^\infty}+\|b^2\|_{L^\infty}\right)\nonumber\\
&\leq& C\left(1+\|F\|_{L^\infty} \right)\leq
C\left(1+\|F_x\|_{L^2}\right)\label{2.28}
\end{eqnarray}
and
\begin{eqnarray}
\|u_{xx}\|_{L^2}&\leq&
C\left(\|F_x\|_{L^2}+\|P(\rho)_x\|_{L^2}+\|(b^2)_x\|_{L^2}\right)\nonumber\\
&\leq& C\left( \|F_x\|_{L^2}
+\|\rho_x\|_{L^2}+\|b_x\|_{L^2}\right).\label{2.29}
\end{eqnarray}

Next, if we differentiate (\ref{1.5})$_1$, (\ref{1.5})$_3$ with
respect to $x$, multiply  the resulting equations by $\rho_x,b_x$ in
$L^2$ respectively, and integrate by parts, then we deduce from
(\ref{2.3}), (\ref{2.23}), (\ref{2.28}) and (\ref{2.29}) that
\begin{eqnarray}
&&\frac{d}{dt}\left(\|\rho_x\|_{L^2}^2+\|b_x\|_{L^2}^2\right)\nonumber\\
&&\quad
\leq C\|u_x\|_{L^\infty} \left(
\|\rho_x\|_{L^2}^2+\|b_x\|_{L^2}^2\right)+C \|u_{xx}\|_{L^2}
\left(\|\rho_x\|_{L^2}+\|b_x\|_{L^2}\right)\nonumber\\
&&\quad\leq C\left(1+\|F_x\|_{L^2}\right)\left(
\|\rho_x\|_{L^2}^2+\|b_x\|_{L^2}^2\right)\nonumber\\
&&\qquad+C\left(\|F_x\|_{L^2}+\|\rho_x\|_{L^2}+\|b_x\|_{L^2}\right)\left(\|\rho_x\|_{L^2}+\|b_x\|_{L^2}\right)\nonumber\\
&&\quad\leq C\left(1+\|F_x\|_{L^2}\right)\left(1+
\|\rho_x\|_{L^2}^2+\|b_x\|_{L^2}^2\right),\label{2.30}
\end{eqnarray}
which, combined with (\ref{2.19}) and Gronwall's inequality, shows
$$
\sup_{0\leq t<
T}\left(\|\rho_x(t)\|_{L^2}^2+\|b_x(t)\|_{L^2}^2\right)\leq C.
$$
Moreover, it readily follows from (\ref{1.5})$_1$ and
(\ref{1.5})$_3$ that $\|\rho_t(t)\|_{L^2},\|b_t(t)\|_{L^2}$ are
bounded on $[0,T)$, since one easily deduces from (\ref{2.2}),
(\ref{2.3}), (\ref{2.18}), (\ref{2.23}) and Sobolev's inequality that
\begin{equation*}
\|u(t)\|_{L^\infty}^2\leq
C\left(\|u(t)\|_{L^2}^2+\|u(t)\|_{L^2}\|u_x(t)\|_{L^2}\right)\leq
C,\quad\forall\ t\in[0, T).
\end{equation*}
The proof of Lemma \ref{lem2.7} is thus finished.\hfill$\square$

\vskip 2mm

Based on Lemmas \ref{lem2.5}--\ref{lem2.7}, we can easily derive higher-order estimates for the velocity.
\begin{lem}\label{lem2.8}Let $(\rho,u,b)$ be a smooth solution of (\ref{1.5}), (\ref{1.6}) on $[0,1]\times[0,T)$. Then,
\begin{equation}
\int_0^T\left(\|u_{xx}\|_{L^2}^2+\|u_x\|_{L^\infty}^4+\|u_t\|_{L^2}^2\right)dt\leq
C,\label{2.31}
\end{equation}
and moreover,
\begin{equation}
\sup_{0\leq t< T}
\left(\|u_{xx}(t)\|_{L^2}^2+\|u_t(t)\|_{L^2}^2\right)+\int_0^T
\|u_{xt}\|_{L^2}^2 dt\leq C.\label{2.32}
\end{equation}
\end{lem}
 \pf First, it follows from (\ref{2.3}), (\ref{2.8}), (\ref{2.9}), (\ref{2.18}),
(\ref{2.19}), (\ref{2.23}) and (\ref{2.27}) that
\begin{eqnarray*}
\int_0^T\left(\|u_{xx}\|_{L^2}^2+\|u_t\|_{L^2}^2\right)dt&\leq&
C\int_0^T\left(\|F_x\|_{L^2}^2+\|P(\rho)_x\|_{L^2}^2+\|(b^2)_x\|_{L^2}^2\right)dt\nonumber\\
&&+C\int_0^T\left(\| \dot u\|_{L^2}^2+\|u\|_{L^\infty}^2\|u_x\|_{L^2}^2\right)\nonumber\\
 &\leq& C+C\int_0^T\left(\|F_x\|_{L^2}^2+\|\rho^{1/2}\dot u\|_{L^2}^2\right)dt\leq
C,
\end{eqnarray*}
which, together with (\ref{2.18}) and Sobolev's inequality, yields
that $\|u_x\|_{L^\infty}^4\in L^1(0,T)$.

To prove (\ref{2.32}), we differentiate (\ref{1.5})$_2$ with respect
to $t$ to get
$$
\rho u_{tt}+\rho uu_{xt}-\lambda
u_{xxt}+\left(P(\rho)+\frac{b^2}{2}\right)_{xt}=-\rho_t(u_t+uu_x)-\rho
u_tu_x,
$$
which, multiplied by $u_t$ in $L^2$ and integrated by parts, results
in
\begin{eqnarray*}
&&\frac{1}{2}\frac{d}{dt}\int_0^1\rho u_t^2dx+\lambda\int_0^1
u_{xt}^2dx\nonumber\\
&&\quad=\int_0^1\left[\left(P(\rho)+\frac{b^2}{2}\right)_{t}u_{xt}-\rho_t(u_t^2+u u_x u_t)-\rho u_t^2u_x\right]dx\nonumber\\
&&\quad \leq
C\int_0^1\left[\left(|\rho_t|+|b_t| \right)|u_{xt}|+|u_t|\left(|\rho_t||u_{t}|+|\rho_t||u_x|+|u_x||u_t|\right)\right]dx\nonumber\\
&&\quad\leq C
\left(\|\rho_t\|_{L^2}+\|b_t\|_{L^2}\right)\|u_{xt}\|_{L^2} +C \|u_t\|_{L^\infty} \left(\|u_t\|_{L^2}^2+\|\rho_t\|_{L^2}^2+\|u_x\|_{L^2}^2\right)\nonumber\\
&&\quad\leq C\|u_{xt}\|_{L^2}+C\left(
\|u_{t}\|_{L^2}+\|u_{xt}\|_{L^2}\right)\left(1+\|u_{t}\|_{L^2}^2\right)\nonumber\\
&&\quad\leq
\frac{\lambda}{2}\|u_{xt}\|_{L^2}^2+C\left(1+\|u_{t}\|_{L^2}^4\right),
\end{eqnarray*}
where the previous lemmas and Cauchy-Schwarz's inequality have been used.
Thus, it follows from (\ref{2.23}), (\ref{2.31}) and Gronwall's inequality that
\begin{equation*}
\sup_{0\leq t< T} \|u_t(t)\|_{L^2}^2 +\int_0^T \|u_{xt}\|_{L^2}^2
dt\leq C.
\end{equation*}
As a consequence, using (\ref{2.3}), (\ref{2.18}), (\ref{2.23}) and
(\ref{2.27}), we infer from (\ref{1.5})$_2$ that
$\|u_{xx}(t)\|_{L^2}$ is bounded on $[0,T)$. The proof of Lemma
\ref{lem2.8} is thus finished.\hfill$\square$

\vskip 2mm

Based on the local existence result and  the global a-priori estimates established in Lemmas
\ref{lem2.1}--\ref{lem2.8}, we can prove Theorem \ref{thm1.1}.

\vskip 2mm

{\it Proof of Theorem \ref{thm1.1}.}  As aforementioned, the
local-in-time solutions can be obtained via the standard fixed point
theorem. Thus, based on the global a-priori estimates established in
Lemmas \ref{lem2.1}--\ref{lem2.8}, we can extend the local solutions
globally in time on $[0,T)$ for any $0<T<\infty$. This proves the
global existence of strong solutions. The uniqueness of strong
solutions can be easily shown by using the standard $L^2$-method,
and the details are omitted here for simplicity. The proof of
Theorem \ref{thm1.1} is therefore complete.\hfill$\square$

\section{The vanishing resistivity limit}
The global existence of strong solutions to (\ref{1.4}),
(\ref{1.10}) stated in the first part of Theorem \ref{thm1.2} can be
established in a manner similar to (indeed much easier) that used
in \cite{CW2,ZX2008} by combining the standard local existence
result and the global a-priori estimates. Thus, in this section we
only focus on the derivation of the uniform bounds stated in
(\ref{1.12}) and (\ref{1.13}), which will suffice in the study of the
vanishing resistivity limit in the second part of Theorem
\ref{thm1.2}. To do this, we assume that $(\rho,u,b)$ is a smooth
solution of (\ref{1.4}), (\ref{1.10}) defined on
$[0,1]\times[0,T)$. For simplicity, we also denote by $C$ a
generic positive constant, which may depend on $\lambda,A,\gamma,T$,
the norms of the initial data $(\rho_0,u_0,b_0)(x)$ and boundary data
$(b_1, b_2)(t)$, but is independent of $\nu$.

First, to deal with the boundary effects, we need the following
formulas of $b_x$ on the boundaries $x=0,1$.
\begin{lem}\label{lem3.1}Let $(\rho,u,b)$ be a smooth solution of (\ref{1.4}), (\ref{1.10}) on $[0,1]\times[0,T)$. Then,
\begin{eqnarray}
\nu b_x(0,t)&=&\nu
\left(b_2(t)-b_1(t)\right)-\partial_t\left[\int_0^1\left(\int_0^x
b(\xi,t)d\xi\right)dx\right]-\int_0^1(u b)(x,t) dx,\label{3.1}\\
\nu b_x(1,t)&=&\nu
\left(b_2(t)-b_1(t)\right)+\partial_t\left[\int_0^1\left(\int_x^1
b(\xi,t)d\xi\right)dx\right]-\int_0^1(u b)(x,t) dx.\label{3.2}
\end{eqnarray}
\end{lem}
\pf  First, thanks to the non-slip boundary condition $u|_{x=0}=0$, we deduce after integrating (\ref{1.4})$_3$ over $(0,x)$ that
\begin{equation}
\nu b_x(0,t)=\nu b_x -\partial_t \int_0^x b(\xi,t)d\xi -ub,\label{3.3}
\end{equation}
and hence, by integrating (\ref{3.3}) with respect to $x$ over $(0,1)$, we obtain
(\ref{3.1}).

Similarly, integrating (\ref{1.4})$_3$ over $(x,1)$  and using the non-slip boundary condition $u|_{x=1}=0$ again,
we find
\begin{equation}
\nu b_x(1,t)=\nu b_x +\partial_t \int_0^x b(\xi,t)d\xi -ub.\label{3.4}
\end{equation}
Thus, an integration of (\ref{3.4}) with respect to $x$ over $(0,1)$ immediately leads to
(\ref{3.2}).\hfill$\square$

\vskip 2mm

In terms of the boundary formulas stated in Lemma \ref{lem3.1}, one
can easily derive the elementary energy estimates similar to the
ones in Lemma \ref{lem2.1}.
\begin{lem}\label{lem3.2}Let $(\rho,u,b)$ be a smooth solution of (\ref{1.4}), (\ref{1.10}) on $[0,1]\times[0,T)$.
Then for any $0\leq t< T$,
\begin{equation}
\int_0^1\rho(x,t)dx=\int_0^1\rho_0(x)dx,\label{3.5}
\end{equation}
and
\begin{equation}
 \int_0^1\left(\frac{1}{2}\rho
u^2+\frac{1}{2}b^2+\frac{A}{\gamma-1}\rho^\gamma\right)(x,t)dx+\int_0^t\left(\lambda\|
u_x\|_{L^2}^2 +\nu\|b_x\|_{L^2}^2\right)ds\leq C.\label{3.6}
\end{equation}
\end{lem}
\pf As a consequence of conservation of mass, one obtains
(\ref{3.5}). To prove (\ref{3.6}), multiplying (\ref{1.4})$_2$,
(\ref{1.4})$_3$ by $u$ and $b$ in $L^2$ respectively,
integrating by parts over $(0,1)\times(0,t)$ with $0\leq t< T$, and using
(\ref{1.4})$_1$ and Lemma \ref{lem3.1}, we deduce
\begin{eqnarray*}
&&\int_0^1\left(\frac{1}{2}\rho
u^2+\frac{1}{2}b^2+\frac{A}{\gamma-1}\rho^\gamma\right)(x,t)dx+\int_0^t\left(\lambda\|
u_x\|_{L^2}^2 +\nu\|b_x\|_{L^2}^2\right)(x,s)ds\nonumber\\
&&\quad =\int_0^1\left(\frac{1}{2}\rho_0
u_0^2+\frac{1}{2}b_0^2+\frac{A}{\gamma-1}\rho_0^\gamma\right)(x)dx+\int_0^t\nu\left(b_2(s)b_x(1,s)-b_1(s)b_x(0,s)\right)ds\nonumber\\
&&\quad\leq C+C\int_0^1|b(x,t)|dx+C\int_0^t\int_0^1\left(
|b_2'||b|+|b_1'| |b|+|u|| b|\right)dx
ds\nonumber\\
&&\quad\leq
C+C\|b(t)\|_{L^2}+C\int_0^t\left(1+\|u\|_{L^\infty}\right)\|b\|_{L^2}d
s\nonumber\\
&&\quad \leq
C+\frac{1}{4}\|b(t)\|_{L^2}^2+C\int_0^t\left(1+\|u_x\|_{L^2}\right)\|b\|_{L^2}d
s\nonumber\\
&&\quad \leq C+\frac{1}{4}\|b(t)\|_{L^2}^2+C\int_0^t
\|b\|_{L^2}^2ds+\frac{\lambda}{2}\int_0^t\|u_x\|_{L^2}^2ds,
\end{eqnarray*}
where we have also used H\"{o}lder's, Sobolev's and Cauchy-Schwarz's
inequalities. This, together with Gronwall's inequality, immediately
gives (\ref{3.6}).\hfill$\square$

 \vskip 2mm

The upper boundedness of the density for the problem (\ref{1.4}), (\ref{1.10})
can be obtained in the same manner as in Lemma \ref{lem2.2}, and we have
\begin{lem}\label{lem3.3}
Let $(\rho,u,b)$ be a smooth solution of (\ref{1.4}), (\ref{1.10})
on $[0,1]\times[0,T)$. Then,
\begin{equation}
0\leq \rho(x,t)\leq C,\quad \forall\
(x,t)\in[0,1]\times[0,T).\label{3.7}
\end{equation}
\end{lem}

Next, we aim to estimate first-order derivatives of the velocity and derive
the higher integrability of the magnetic field stated in
Lemma \ref{lem2.5}.  We remark that it is more complicated to achieve this due to the
additional boundary effects induced by the magnetic diffusion term $\nu b_{xx}$.
So, to be proceeded, as in Section 2, let `` $\dot{}$ " be the material derivative $\dot f= f_t+uf_x$ and define
\begin{equation}
F\triangleq\left(\lambda
u_x-P(\rho)-\frac{1}{2}b^2\right)(x,t)\quad{\rm with}\quad
F_x=\rho\dot u.\label{3.8}
\end{equation}
where $(\rho,u,b)$ is the solution of (\ref{1.4}), (\ref{1.10}).

Based on a full use of the `` effective viscous flux" $F$ and subtle
boundary analysis, we can show the following estimate in a way similar  to
the proof of Lemma \ref{lem2.5}
\begin{lem}\label{lem3.4}
Let $(\rho,u,b)$ be a smooth solution of (\ref{1.4}), (\ref{1.10})
on $[0,1]\times[0,T)$. Then,
\begin{eqnarray}
&&\sup_{0\leq t<
T}\left(\|u_x(t)\|_{L^2}^2+\|b(t)\|_{L^4}^4+\nu\|b_x(t)\|_{L^2}^2\right)\nonumber\\
&&\qquad+\int_0^T\left(\|\rho^{1/2}\dot
u\|_{L^2}^2+\|b\|_{L^6}^6+\nu\|bb_x\|_{L^2}^2+\nu^2\|b_{xx}\|_{L^2}^2\right)dt\leq
C.\label{3.9}
\end{eqnarray}
\end{lem}
\pf First, similarly to the derivation of (\ref{2.11}), (\ref{2.13})
and (\ref{2.15}),  we deduce from (\ref{1.4})$_1$ and (\ref{1.4})$_3$ that
\begin{eqnarray}
&&\frac{\lambda}{2}\frac{d}{dt}\int_0^1u_x^2dx+\int_0^1\rho \dot u^2 dx  \nonumber\\
&&\quad=-\int_0^1 P(\rho)_x(u_t+uu_x) dx-\frac{1}{2}\int_0^1
(b^2)_x(u_t+uu_x)dx -\frac{\lambda}{2}\int_0^1 u_x^3
dx\nonumber\\
&&\quad=\frac{d}{dt}\int_0^1 \left(P(\rho) + \frac{b^2}{2} \right)
u_{x} dx+\int_0^1\left(\gamma P(\rho) u_x^2 +b^2u_x^2-\nu b b_{xx}
u_x-\frac{\lambda}{2}
u_x^3  \right)dx\nonumber\\
&&\quad \leq \frac{d}{dt}\int_0^1 \left(P(\rho) +
\frac{b^2}{2}\right) u_{x} dx+\frac{\nu^2}{2}\int_0^1 b_{xx}^2 dx+
C\int_0^1\left(
|u_x|^3+b^6\right)dx+C\nonumber\\
&&\quad \leq \frac{d}{dt}\int_0^1 \left(P(\rho)+
\frac{b^2}{2}\right) u_{x} dx+\frac{\nu^2}{2}\int_0^1 b_{xx}^2dx +
C\int_0^1\left( |F|^3+b^6\right)dx+C,\label{3.10}
\end{eqnarray}
where we have also used (\ref{3.7}), (\ref{3.8}) and Cauchy-Schwarz's inequality.

Secondly, analogously to the proof of (\ref{2.17}), we have
\begin{eqnarray*}
&&\frac{d}{dt}\int_0^1 b^4 dx+ \frac{3}{2\lambda}\int_0^1 b^6 dx
+12\nu\int_0^1 b^2b_x^2dx\nonumber\\
&&\quad =
-\frac{3}{\lambda}\int_0^1\left(F+P(\rho) \right) b^4dx+4\nu \left[b_2^3(t)b_x(1,t)-b_1^3(t)b_x(0,t)\right] \nonumber\\
&&\quad \leq C+\frac{1}{2\lambda}\int_0^1 b^6 dx+C\int_0^1 |F|^3
dx+R_1(t),
\end{eqnarray*}
and consequently,
\begin{equation}
\frac{d}{dt}\int_0^1 b^4 dx+  \frac{1}{\lambda}\int_0^1 b^6 dx
+12\nu\int_0^1b^2b_x^2dx\leq C+C\int_0^1 |F|^3
dx+R_1(t),\label{3.11}
\end{equation}
where $R_1(t)$ denotes the boundary term:
$$
R_1(t)\triangleq 4\nu
\left[b_2^3(t)b_x(1,t)-b_1^3(t)b_x(0,t)\right].
$$

Moreover, multiplying (\ref{1.4})$_3$ by $2\nu b_{xx}$ in $L^2$ and
integrating by parts, we use the non-slip boundary conditions $u|_{x=0,1}=0$ to deduce
\begin{eqnarray}
&&\nu\frac{d}{dt}\int_0^1 b_x^2dx+2\nu^2\int_0^1 b_{xx}^2dx\nonumber\\
&&\quad =2\nu\int_0^1\left(u_xb+ub_x
\right)b_{xx}dx+2\nu\left[b_2'(t)
b_x(1,t)-b_1'(t)b_x(0,t)\right]\nonumber\\
&&\quad =\nu\int_0^1\left(2u_xb b_{xx}-u_xb_x^2
\right)dx+2\nu\left[b_2'(t)
b_x(1,t)-b_1'(t)b_x(0,t)\right].\label{3.12}
\end{eqnarray}

Due to non-negativity of $P(\rho)$ and $b^2$, it follows from (\ref{3.8}) that
\begin{equation*}
-\int_0^1 u_x
b_x^2dx=-\frac{1}{\lambda}\int_0^1\left(F+P(\rho)+\frac{b^2}{2}\right)b_x^2dx
\leq -\frac{1}{\lambda}\int_0^1 F b_x^2dx,
\end{equation*}
which, inserted into (\ref{3.12}) and combined with (\ref{3.7}),
(\ref{3.8}) and Cauchy-Schwarz's inequality, yields
\begin{eqnarray*}
&&\nu\frac{d}{dt}\int_0^1 b_x^2dx+2\nu^2\int_0^1 b_{xx}^2dx
\nonumber\\
&&\quad \leq 2\nu\int_0^1 u_xb b_{xx}dx-\frac{\nu}{\lambda}\int_0^1 F
b_x^2dx+R_2(t)\nonumber\\
&& \quad\leq
\frac{\nu^2}{2}\int_0^1b_{xx}^2dx+C\int_0^1\left(u_x^3+b^6+\nu |F|b_x^2 \right)dx+R_2(t)\nonumber\\
&&\quad \leq
\frac{\nu^2}{2}\int_0^1b_{xx}^2dx+C+C\int_0^1\left(F^3+b^6+\nu|F|b_x^2\right)
dx+R_2(t),
\end{eqnarray*}
where $R_2(t)$ is the boundary term:
$$
R_2(t)\triangleq 2\nu\left[b_2'(t) b_x(1,t)-b_1'(t)b_x(0,t)\right].
$$

Thus, we obtain
\begin{equation}
\nu\frac{d}{dt}\int_0^1 b_x^2dx+\frac{3\nu^2}{2}\int_0^1
b_{xx}^2dx\leq C+C\int_0^1\left(F^3+b^6+\nu|F|b_x^2\right)
dx+R_2(t).\label{3.13}
\end{equation}
Adding (\ref{3.10}) and (\ref{3.13}) together, we arrive at
\begin{eqnarray}
&&\frac{d}{dt}\int_0^1\left(\frac{\lambda}{2}u_x^2+\nu
b_x^2\right)dx+\int_0^1\left(\rho\dot u^2+\nu^2
b_{xx}^2\right)dx\nonumber\\
&&\quad \leq \frac{d}{dt}\int_0^1 \left(P(\rho) + \frac{b^2}{2}
\right) u_{x} dx+C\int_0^1\left(F^3+b^6+\nu |F|b_x^2\right)
dx+R_2(t)+C.\label{3.14}
\end{eqnarray}

Note that it follows from Sobolev's and Young's inequalities that for
any $0\leq t< T$,
\begin{eqnarray}
|R_1(t)|+|R_2(t)|&\leq&C\nu\|b_x\|_{L^\infty}\leq
C\nu\left(\|b_x\|_{L^2}+\|b_x\|_{L^2}^{1/2}\|b_{xx}\|_{L^2}^{1/2}\right)\nonumber\\
&\leq&\varepsilon\nu^{2}\|b_{xx}\|_{L^2}^2+C\varepsilon^{-1}\left(\nu
\|b_x\|_{L^2}^2 +\nu^{1/2}\right),\label{3.15}
\end{eqnarray}
where $C>0$ depends on the $C^1$-norm of $b_1(t)$ and $b_2(t)$, but
is independent of $\varepsilon,\nu\in(0,1)$.

Now, multiplying (\ref{3.11}) a sufficiently large number $K>0$, adding
the resulting inequality to (\ref{3.14}) and integrating it over
$(0,t)$, by (\ref{3.15}) with $0<\varepsilon<1$ being sufficiently
small, we deduce in a manner similar to that used in the derivation of (\ref{2.20}) that for $\forall$ $\nu\in(0,1)$,
\begin{eqnarray}
&&\left(\|u_x(t)\|_{L^2}^2+\|b(t)\|_{L^4}^4+\nu\|b_x(t)\|_{L^2}^2\right)\nonumber\\
&&\qquad+\int_0^t\left(\|\rho^{1/2}\dot
u\|_{L^2}^2+\|b\|_{L^6}^6+\nu\|bb_x\|_{L^2}^2+\nu^2\|b_{xx}\|_{L^2}^2\right)ds\nonumber\\
&&\quad\leq
C+C\int_0^t\|F\|_{L^3}^3ds+C\nu\int_0^t\|F\|_{L^\infty}\|b_x\|_{L^2}^2ds\nonumber\\
&&\quad\leq
C+C \nu^2\int_0^t\|b_x\|_{L^2}^4ds+C\int_0^t\left(\|F\|_{L^3}^3+\|F\|_{L^\infty}^2\right)ds\nonumber\\
&&\quad\leq
C+C \nu^2\int_0^t\|b_x\|_{L^2}^4ds+\frac{1}{2}\int_0^t\|\rho^{1/2}\dot u\|_{L^2}^2ds\nonumber\\
&&\qquad+C \int_0^t(1+\|u_x\|_{L^2}^2)\left(1
+\|u_x\|_{L^2}^2+\|b\|_{L^4}^4\right)  ds,\label{3.16}
\end{eqnarray}
where (\ref{3.6})--(\ref{3.8}), Sobolev's, H\"{o}lder's and
Cauchy-Schwarz's inequalities have been used to get
\begin{eqnarray*}
\|F\|_{L^3}^3+\|F\|_{L^\infty}^2 &\leq&
C\left(1+\|F\|_{L^1}\right)\left(\|F\|_{L^2}^2+\|F\|_{L^2}\|F_x\|_{L^2}\right)
\nonumber\\
&\leq&
C\left(1+\|u_x\|_{L^2}\right)\left(1+\|u_x\|_{L^2}^2+\|b\|_{L^4}^4\right)\nonumber\\
&&+C\left(1+\|u_x\|_{L^2}\right)\left(1+\|u_x\|_{L^2}+\|b\|_{L^4}^2\right)\|\rho^{1/2}\dot
u\|_{L^2}
\nonumber\\
&\leq&\frac{1}{2}\|\rho^{1/2} \dot u\|_{L^2}^2+C
\left(1+\|u_x\|_{L^2}^2\right)\left(1
+\|u_x\|_{L^2}^2+\|b\|_{L^4}^4\right),
\end{eqnarray*}
which is analogous to (\ref{2.22}). As an immediate result of (\ref{3.16}), we have
\begin{eqnarray*}
&&\left(\|u_x(t)\|_{L^2}^2+\|b(t)\|_{L^4}^4+\nu\|b_x(t)\|_{L^2}^2\right)\nonumber\\
&&\qquad+\int_0^t\left(\|\rho^{1/2}\dot
u\|_{L^2}^2+\|b\|_{L^6}^6+\nu\|bb_x\|_{L^2}^2+\nu^2\|b_{xx}\|_{L^2}^2\right)ds\nonumber\\
&&\quad\leq C+C \nu^2\int_0^t\|b_x\|_{L^2}^4ds +C
\int_0^t(1+\|u_x\|_{L^2}^2)\left(1
+\|u_x\|_{L^2}^2+\|b\|_{L^4}^4\right) ds.
\end{eqnarray*}
The above inequality, combined with Gronwall's inequality and the fact that
$\|u_x\|_{L^2}^2+\nu\|b_x\|_{L^2}^2\in L^1(0,T)$ due to (\ref{3.6}),
leads to the desired estimate (\ref{3.9}).\hfill$\square$

\vskip 2mm

By virtue of (\ref{3.9}), the upper boundedness of the magnetic field and the
lower boundedness of the density can be obtained in a manner similar to that used in the proof of (\ref{2.23}).
\begin{lem}\label{lem3.5}
Let $(\rho,u,b)$ be a smooth solution of (\ref{1.4}), (\ref{1.10}) on $[0,1]\times[0,T)$. Then,
\begin{equation}
\sup_{0\leq t<
T}\left(\|b(t)\|_{L^\infty}+\|\rho^{-1}(t)\|_{L^\infty}\right)\leq
C.\label{3.17}
\end{equation}
\end{lem}
\pf Similarly to the proof of (\ref{2.24}), we multiply (\ref{1.4})$_3$ by $2n b^{2n-1}$ with $1\leq n\in \mathbb{N}$ and
integrate by parts over $(0,1)$ to infer, keeping in mind the non-negativity of
$P(\rho)$ and $b^2$, that there exists a positive constant $C$, independent of $n$ and $\nu$, such that
\begin{eqnarray*}
&&\frac{d}{dt}\int_0^1b^{2n}(x,t)dx+2n(2n-1)\nu\int_0^1
b^{2n-2}b_x^2dx\nonumber\\
&&\quad=-(2n-1)\int_0^1 u_x
b^{2n}dx+2n\nu\left[b_2^{2n-1}(t)b_x(1,t)-b_1^{2n-1}(t)b_x(0,t)\right]\nonumber\\
&&\quad\leq -\frac{2n-1}{\lambda}\int_0^1
\left(F+P(\rho)+\frac{b^2}{2}\right)
b^{2n}dx+C n C^{2n}\nu\|b_x\|_{L^\infty}\nonumber\\
&&\quad\leq -\frac{2n-1}{\lambda}\int_0^1 F
b^{2n}dx+C n C^{2n}\nu\|b_x\|_{L^\infty}\nonumber\\
&&\quad\leq C n\|F\|_{L^\infty}\|b\|_{L^{2n}}^{2n}
+C n C^{2n}\nu\left(\|b_x\|_{L^2}+\|b_x\|_{L^2}^{1/2}\|b_{xx}\|_{L^2}^{1/2}\right)\nonumber\\
&&\quad\leq C n\|F\|_{L^\infty}\|b\|_{L^{2n}}^{2n} +C n
C^{2n}\left(1+\nu^2\|b_{xx}\|_{L^2}^2\right),
\end{eqnarray*}
where we have used (\ref{3.9}), Sobolev's and Cauchy-Schwarz's
inequalities. Noting that (\ref{3.7}), together with (\ref{3.8}) and (\ref{3.9}), gives
\begin{eqnarray}
\|F\|_{L^\infty}^2&\leq& \|F\|_{L^2}^2+\|F\|_{L^2}\|F_x\|_{L^2}\leq
C\left(1+\|F_x\|_{L^2}\right)\nonumber\\
&\leq& C\left(1+\|\rho^{1/2}\dot u\|_{L^2}\right)\in L^2(0,T),\label{3.18}
\end{eqnarray}
thus we have by (\ref{3.9}) that
\begin{eqnarray*}
\|b(t)\|_{L^{2n}}^{2n}&\leq& C n
C^{2n}\exp\left\{Cn\int_0^T\|F\|_{L^\infty}dt\right\}\int_0^T\left(1+\nu^2\|b_{xx}\|_{L^2}^2dt\right)\nonumber\\
&\leq& C nC^{2n}\exp\left\{Cn\right\},
\end{eqnarray*}
from which we conclude that $\|b(t)\|_{L^\infty}\leq C$ for all $t\in[0, T)$ by raising to the power $1/(2n)$
to both sides and letting $n\to\infty$. With this at hand, the lower boundedness of the density
can be shown exactly in the same way as in Lemma \ref{lem2.6}. The proof of (\ref{3.17}) is complete.
\hfill$\square$

\vskip 2mm

The following refined estimates of density and magnetic field play
an important role in the analysis of vanishing resistivity limit,
and will be also used to deal with the thickness of magnetic boundary layers.
\begin{lem}
\label{lem3.6}For any $ 0<\nu<1$, there exists a positive constant
$C$, independent of $\nu$, such that
\begin{equation}
\nu^{1/2}\sup_{0\leq t<
T}\left(\|b_x(t)\|_{L^2}^2+\|\rho_x(t)\|_{L^2}^2\right)
+\int_0^T\left(\nu^{3/2}\| b_{xx}\|_{L^2}^2+\nu^{1/2}\|b
b_x\|_{L^2}^2\right)dt \leq C.\label{3.19}
\end{equation}
\end{lem}
\pf Differentiating (\ref{1.4}) with respect to $x$, multiplying the
resulting equation by $2\rho_x$ in $L^2$ and integrating by parts,
by (\ref{3.7}), (\ref{3.8}), (\ref{3.17}) and (\ref{3.18}) we deduce
\begin{eqnarray}
\frac{d}{dt}\int_0^1 \rho_x^2 dx&=&-3\int_0^1 u_x \rho_x^2 dx-2
\int_0^1
\rho \rho_xu_{xx} dx\nonumber\\
&=&-3\int_0^1 u_x \rho_x^2 dx-\frac{2}{\lambda}\int_0^1
\left(F_x+P_x+ b b_x\right)\rho \rho_x dx\nonumber\\
&\leq&C\left(1+\|u_x\|_{L^\infty}\right) \|\rho_x\|_{L^2}^2 +
C\left(\|F_x\|_{L^2}+\|\rho_x\|_{L^2}+\|b_x\|_{L^2}\right)\|\rho_x\|_{L^2}\nonumber\\
&\leq&C\left(1
+\|F_x\|_{L^2}\right)\left(\|\rho_x\|_{L^2}^2+\|b_x\|_{L^2}^2\right)+C\left(1
+\|F_x\|_{L^2}^2\right),\label{3.20}
\end{eqnarray}
where we have used the following estimate which follows from (\ref{3.7}), (\ref{3.8}), (\ref{3.17}) and
(\ref{3.18}).
\begin{equation}
\|u_x\|_{L^\infty}\leq
C\left(\|F\|_{L^\infty}+\|P(\rho)\|_{L^\infty}+\|b\|_{L^\infty}^2\right)
\leq C\left(1+\|F_x\|_{L^2}\right).\label{3.21}
\end{equation}

Using the non-slip boundary conditions and integrating by parts, we get from (\ref{3.12}) that
\begin{eqnarray}
&&\nu\frac{d}{dt}\int_0^1 b_x^2dx+2\nu^2\int_0^1 b_{xx}^2dx =2\nu\int_0^1\left(u_xb+ub_x \right)b_{xx}dx+R_2(t)\nonumber\\
&&\quad =-2\nu\int_0^1 u_{xx}b b_{x}dx-3\nu\int_0^1u_x b_x^2
dx+R_2(t)+R_3(t),\label{3.22}
\end{eqnarray}
where $R_2(t)$ is the same as in (\ref{3.13}) and $R_3(t)$ is also a boundary term denoted by
$$
R_3(t)\triangleq 2\nu\left[
b_2(t)(u_xb_x)(1,t)-b_1(t)(u_xb_x)(0,t)\right].
$$

We are now in a position to bound the terms on the right-hand side of (\ref{3.22}). First, recalling the definition of $F$ in
(\ref{3.8}), using (\ref{3.7}) and (\ref{3.17}), we have
\begin{eqnarray}
-2\nu\int_0^1 u_{xx}b b_{x}dx&=&-\frac{2\nu}{\lambda}\int_0^1 \left(F_{x}+P(\rho)_x+bb_x\right)b b_{x}dx\nonumber\\
&\leq & C\nu\left(\|F_x\|_{L^2}+\|\rho_x\|_{L^2}+\|b_x\|_{L^2}\right)\| b_x\|_{L^2}   \nonumber\\
&\leq & C\nu\left(\|F_x\|_{L^2}^2+\|\rho_x\|_{L^2}^2+\|b_x\|_{L^2}^2\right),\label{3.23}
\end{eqnarray}
and by (\ref{3.21}) we find
\begin{equation}
-3\nu\int_0^1u_x b_x^2 dx\leq
C\nu\|u_x\|_{L^\infty}\|b_x\|_{L^2}^2\leq
C\nu\left(1+\|F_x\|_{L^2}\right)\|b_x\|_{L^2}^2.\label{3.24}
\end{equation}

Choosing $0<\varepsilon<1$ small enough, we infer from
(\ref{3.15}) that
\begin{equation}
R_2(t)\leq \frac{\nu^2}{2}\|b_{xx}\|_{L^2}^2+C\nu \|b_x\|_{L^2}^2+C
\nu^{1/2}.\label{3.25}
\end{equation}

In order to deal with $R_3(t)$,  we observe with the help of (\ref{3.21}) that
\begin{equation*}
|u_x(0,t)|+|u_x(1,t)|\leq C\|u_x\|_{L^\infty}\leq
C\left(1+\|F_x\|_{L^2}\right),
\end{equation*}
and hence,
\begin{eqnarray}
R_3(t)&\leq& C\nu\left(1+\|F_x\|_{L^2}\right)\left(\|b_1\|_{C([0,T])}+\|b_2\|_{C([0,T])}\right)\|b_x\|_{L^\infty}\nonumber\\
&\leq&C\nu\left(1+\|F_x\|_{L^2}\right)\left(\|b_x\|_{L^2}+\|b_x\|_{L^2}^{1/2}\|b_{xx}\|_{L^2}^{1/2}\right)\nonumber\\
&\leq&\frac{\nu^2}{2}\|b_{xx}\|_{L^2}^2+C\nu
\|b_x\|_{L^2}^2+C\left(\nu+\nu^{1/2}\right)\left(1+\|F_x\|_{L^2}^2\right).\label{3.26}
\end{eqnarray}

Substituting (\ref{3.23})--(\ref{3.26}) into (\ref{3.22}), we conclude that for  any $\nu\in(0,1)$,
\begin{eqnarray}
&&\nu\frac{d}{dt}\|b_x\|_{L^2}^2+\left(\nu^2\| b_{xx}\|_{L^2}^2+\nu\|b b_x\|_{L^2}^2\right)\nonumber\\
&&\quad \leq C\nu\left(\|b_x\|_{L^2}^2+\|\rho_x\|_{L^2}^2\right)\left(1+\|F_x\|_{L^2}\right) +C
\nu^{1/2} \left(1+\|F_x\|_{L^2}^2\right).  \label{3.27}
\end{eqnarray}

Now, multiplying (\ref{3.20}) by $\nu$ and adding the resulting
inequality to (\ref{3.27}), we obtain that for any $0<\nu<1$,
\begin{eqnarray*}
&&\nu\frac{d}{dt}\left(\|\rho_x\|_{L^2}^2+\|b_x\|_{L^2}^2\right)+ \nu^2\| b_{xx}\|_{L^2}^2 \nonumber\\
&&\quad \leq
C\nu\left(\|\rho_x\|_{L^2}^2+\|b_x\|_{L^2}^2\right)\left(1+\|F_x\|_{L^2}\right)+C\nu^{1/2}
\left(1+\|F_x\|_{L^2}^2\right),
\end{eqnarray*}
which, combined with (\ref{3.9}), (\ref{3.18}) and Gronwall's inequality, implies
\begin{eqnarray*}
&&\nu\sup_{0\leq t\leq T} \left(\|b_x\|_{L^2}^2+\|P(\rho)_x\|_{L^2}^2\right)+\int_0^T (\nu^2\|
b_{xx}\|_{L^2}^2+\nu\|b b_x\|_{L^2}^2) dt \nonumber\\
&& \quad \leq C\nu^{1/2}\exp\left\{\int_0^T\left(1+\|F_x\|_{L^2}\right)dt\right\}
\int_0^T\left(1+\|F_x\|_{L^2}^2\right)dt\\
&&\quad \leq C\nu^{1/2}\exp\left\{\int_0^T(1+\|\rho^{1/2}\dot u\|_{L^2})dt\right\}
\int_0^T (1+\|\rho^{1/2}\dot u\|_{L^2}^2) dt\\
&&\quad \leq C\nu^{1/2},
\end{eqnarray*}
which immediately leads to  the desired estimate (\ref{3.19}). \hfill$\square$

\vskip 2mm

With the estimates established in Lemmas \ref{lem2.1}--\ref{lem2.8}
and \ref{lem3.2}--\ref{lem3.6} at hand, we are now ready to prove the second part of Theorem \ref{thm1.2}.

\vskip 2mm

{\it Proof of (ii)-Theorem \ref{thm1.2}.} Let
$(\rho^\nu,u^\nu,b^\nu)$ and $(\rho,u,b)$ be the solutions to the
initial-boundary value problems (\ref{1.4}), (\ref{1.10}) and
(\ref{1.5}), (\ref{1.6}), respectively. Then, by a straightforward
calculation we find that $(\rho^\nu-\rho,u^\nu-u,b^\nu-b)$ satisfies
\begin{equation}
\left(\rho^\nu-\rho\right)_t+\left(\rho^\nu-\rho\right)u^\nu_x+\rho\left(u^\nu-u\right)_x
+\left(\rho^\nu-\rho\right)_xu^\nu+\rho_x\left(u^\nu-u\right)=0,\label{3.28}
\end{equation}
\begin{eqnarray}
\rho^\nu\left(u^\nu-u\right)_t&+&\rho^\nu
u^\nu\left(u^\nu-u\right)_x-\lambda\left(u^\nu-u\right)_{xx}
=-\left(\rho^\nu-\rho\right)u_t-\left(\rho^\nu-\rho\right)u
u_x\nonumber\\
&-&\rho^\nu\left(u^\nu-u\right)u_x-\left[P(\rho^\nu)-P(\rho)\right]_x-\frac{1}{2}\left[(b^\nu)^2-b^2\right]_x,\label{3.29}
\end{eqnarray}
and
\begin{equation}
\left(b^\nu-b\right)_t+u^\nu_x\left(b^\nu-b\right)+b\left(u^\nu-u\right)_x+u^\nu\left(b^\nu-b\right)_x+\left(u^\nu-u\right)
b_x =\nu  b^\nu_{xx}.\label{3.30}
\end{equation}

First, multiplying (\ref{3.30}) by $2(b^\nu-b)$ in $L^2$ and
integrating by parts, we obtain
\begin{eqnarray}
\frac{d}{dt}\|b^\nu-b\|_{L^2}^2&=&-\int_0^1u^\nu_x\left(b^\nu-b\right)^2dx-2\int_0^1b\left(u^\nu-u\right)_x\left(b^\nu-b\right)dx\nonumber\\
&&-2\int_0^1\left(u^\nu-u\right) b_x\left(b^\nu-b\right)dx+2\nu\int_0^1 b^\nu_{xx}\left(b^\nu-b\right)dx\nonumber\\
&\triangleq&I_1+I_2+I_3+I_4.\label{3.31}
\end{eqnarray}

On one hand, using (\ref{2.23}), (\ref{2.27}), (\ref{3.18}),
(\ref{3.21}) and Sobolev's inequality, we can control the first three
terms on the right-hand side of (\ref{3.31}) as follows.
\begin{eqnarray}
\sum_{i=1}^3I_i&\leq&
C\|u^\nu_x\|_{L^\infty}\|b^\nu-b\|_{L^2}^2+C \|b\|_{L^\infty}\|(u^\nu-u)_x\|_{L^2}\|b^\nu-b\|_{L^2}\nonumber\\
&&+\|b_x\|_{L^2}\|u^\nu-u\|_{L^\infty}\|b^\nu-b\|_{L^2}\nonumber\\
&\leq& C\left(1+\|\sqrt{\rho^\nu}\dot
u^\nu\|_{L^2}\right)\|b^\nu-b\|_{L^2}^2+C\|(u^\nu-u)_x\|_{L^2}
\|b^\nu-b\|_{L^2}\nonumber\\
&\leq&
\varepsilon\|(u^\nu-u)_x\|_{L^2}^2+C\varepsilon^{-1}\left(1+\|\sqrt{\rho^\nu}\dot
u^\nu\|_{L^2}\right)\|b^\nu-b\|_{L^2}^2.\label{3.32}
\end{eqnarray}
On the other hand, by (\ref{2.23}), (\ref{2.27}) and Young's
inequality we deduce
\begin{eqnarray}
I_4&=&-2\nu\int_0^1 b^\nu_x\left(b^\nu_x-b_x\right)  dx+2\nu \left[ b_x^\nu(b^\nu-b)(1,t)-b_x^\nu(b^\nu-b)(0,t)\right]\nonumber\\
&\leq& -2\nu\|b^\nu_x\|_{L^2}^2+2\nu\|b^\nu_x\|_{L^2}\|b_x\|_{L^2}+C\nu\left(\|b\|_{L^\infty}+\|(b_1,b_2)\|_{C([0,T])}\right)\|b^\nu_x\|_{L^\infty}\nonumber\\
&\leq&-\nu\|b^\nu_x\|_{L^2}^2+C\nu+C\nu\left(\|b^\nu_x\|_{L^2}+\|b^\nu_x\|_{L^2}^{1/2}\|b^\nu_{xx}\|_{L^2}^{1/2}\right)\nonumber\\
&\leq&
-\frac{\nu}{2}\|b^\nu_x\|_{L^2}^2+C\nu^2\|b^\nu_{xx}\|_{L^2}^2+C\nu^{1/2}.\label{3.33}
\end{eqnarray}

Putting (\ref{3.32}) and (\ref{3.33}) into (\ref{3.31}), we obtain
($0<\varepsilon<1$)
\begin{eqnarray}
\frac{d}{dt}\|b^\nu-b\|_{L^2}^2+\nu\|b^\nu_x\|_{L^2}^2&\leq&
C\left(\nu^{1/2}+\nu^2\|b^\nu_{xx}\|_{L^2}^{2}\right)+
\varepsilon\|(u^\nu-u)_x\|_{L^2}^2\nonumber\\
&&+C\varepsilon^{-1}\left(1+ \|\sqrt{\rho^\nu}\dot u^\nu\|_{L^2}
\right)\|b^\nu-b\|_{L^2}^2.\label{3.34}
\end{eqnarray}

Similarly, multiplying (\ref{3.28})  by $2(\rho^\nu-\rho)$ in $L^2$,
integrating by parts, and using Lemmas \ref{2.2}--\ref{lem2.8},
(\ref{3.18}), (\ref{3.21}) and Sobolev's inequality, we deduce
\begin{eqnarray}
\frac{d}{dt}\|\rho^\nu-\rho\|_{L^2}^2&\leq&
C\|u^\nu_x\|_{L^\infty}\|\rho^\nu-\rho\|_{L^2}^2+C \|\rho\|_{L^\infty}\|(u^\nu-u)_x\|_{L^2}\|\rho^\nu-\rho\|_{L^2}\nonumber\\
&&+\|\rho_x\|_{L^2}\|u^\nu-u\|_{L^\infty}\|\rho^\nu-\rho\|_{L^2}\nonumber\\
&\leq& C\left(1+\|\sqrt{\rho^\nu}\dot
u^\nu\|_{L^2}\right)\|\rho^\nu-\rho\|_{L^2}^2+C\|(u^\nu-u)_x\|_{L^2}
\|\rho^\nu-\rho\|_{L^2}\nonumber\\
&\leq&
\varepsilon\|(u^\nu-u)_x\|_{L^2}^2+C\varepsilon^{-1}\left(1+\|\sqrt{\rho^\nu}\dot
u^\nu\|_{L^2}\right)\|\rho^\nu-\rho\|_{L^2}^2.\label{3.35}
\end{eqnarray}

Finally, multiplying (\ref{3.29})  by $2(u^\nu-u)$ in $L^2$ and
integrating by parts, by virtue of Lemmas \ref{2.2}--\ref{lem2.8},
(\ref{3.7}), (\ref{3.17}) and Sobolev's inequality again, we see that
\begin{eqnarray}
&&\frac{d}{dt}\|\sqrt{\rho^\nu}(u^\nu-u)\|_{L^2}^2+2\lambda\|(u^\nu-u)_x\|_{L^2}^2\nonumber\\
&&\quad\leq
C \left(\|u_t\|_{L^2}+\|u\|_{L^\infty}\|u_x\|_{L^2}\right)\|\rho^\nu-\rho\|_{L^2}\|u^\nu-u\|_{L^\infty}
+C\|u_x\|_{L^\infty}\|\sqrt{\rho^\nu}(u^\nu-u)\|_{L^2}^2\nonumber\\
&&\qquad+C \left(\|P(\rho^\nu)-P(\rho)\|_{L^2}+ \|(b^\nu)^2-b^2\|_{L^2}\right) \|(u^\nu-u)_x\|_{L^2}\nonumber\\
&&\quad\leq C\left(\|\rho^\nu-\rho\|_{L^2}+ \|b^\nu-b\|_{L^2}
\right)\|(u^\nu-u)_x\|_{L^2}+C\|u_x\|_{L^\infty}\|\sqrt{\rho^\nu}(u^\nu-u)\|_{L^2}^2\nonumber\\
&&\quad\leq \varepsilon\|(u^\nu-u)_x\|_{L^2}^2+
C\varepsilon^{-1}\left(\|\rho^\nu-\rho\|_{L^2}^2+
\|b^\nu-b\|_{L^2}^2
\right)\nonumber\\
&&\qquad+C\left(\|u_x\|_{L^2}+\|u_{xx}\|_{L^2}\right)\|\sqrt{\rho^\nu}(u^\nu-u)\|_{L^2}^2\nonumber\\
&&\quad\leq \varepsilon\|(u^\nu-u)_x\|_{L^2}^2+
C\varepsilon^{-1}\left(\|\rho^\nu-\rho\|_{L^2}^2+
\|b^\nu-b\|_{L^2}^2
\right)+C\|\sqrt{\rho^\nu}(u^\nu-u)\|_{L^2}^2.\label{3.36}
\end{eqnarray}

Thus, putting (\ref{3.34}), (\ref{3.35}) and (\ref{3.36}) together, taking $\varepsilon>0$ suitably small,
we arrive at
\begin{eqnarray*}
&&\frac{d}{dt}\left(\|\rho^\nu-\rho\|_{L^2}^2+
\|b^\nu-b\|_{L^2}^2+\|\sqrt{\rho^\nu}(u^\nu-u)\|_{L^2}^2\right)+\|(u^\nu-u)_x\|_{L^2}^2+\nu\|b^\nu_x\|_{L^2}^2\nonumber\\
&&\quad\leq C\left(1+\|\sqrt{\rho^\nu}\dot u^\nu\|_{L^2}
\right)\left(\|\rho^\nu-\rho\|_{L^2}^2+
\|b^\nu-b\|_{L^2}^2+\|\sqrt{\rho^\nu}(u^\nu-u)\|_{L^2}^2\right)\nonumber\\
&&\qquad +C\left(\nu^{1/2}+\nu^2\|b^\nu_{xx}\|_{L^2}^{2}\right),
\end{eqnarray*}
which, combined with (\ref{3.9}), (\ref{3.19}) and Gronwall's inequality, yields
\begin{eqnarray*}
&&\sup_{0\leq t< T}\left(\|\rho^\nu-\rho\|_{L^2}^2+
\|b^\nu-b\|_{L^2}^2+\|\sqrt{\rho^\nu}(u^\nu-u)\|_{L^2}^2\right)(t)   \nonumber\\
&&\qquad+\int_0^T\left(\|(u^\nu-u)_x\|_{L^2}^2+\nu\|b^\nu_x\|_{L^2}^2\right)dt\leq C\nu^{1/2}.
\end{eqnarray*}
This inequality, together with the
strictly positive lower bound of density due to (\ref{3.17}), justifies the
vanishing resistivity limit stated in the second part of Theorem \ref{thm1.2}.\hfill$\square$

\section{The thickness of magnetic boundary layers}
This section is concerned with the thickness of magnetic boundary layers.
As aforementioned in Section 1, to simplify the analysis, we assume that the one-dimensional
compressible isentropic viscous and non-resistive MHD system
(\ref{1.5}) is equipped with the specific initial and boundary data as follows:
\begin{equation}
(\rho,u,b)|_{t=0}=(\overline\rho,0,0)\quad{\rm and}\quad
u|_{x=0,1}=0,\label{4.1}
\end{equation}
where $\overline\rho\equiv{\rm const.}$ is a positive constant. Then
thanks to the uniqueness result, we conclude that the
initial-boundary value problem (\ref{1.5}), (\ref{4.1}) has only a
trivial solution:
\begin{equation}
(\rho,u,b)(x,t)=(\overline\rho,0,0),\quad(x,t)\in[0,1]\times[0,T).\label{4.2}
\end{equation}
This particularly finishes the proof of the first part of Theorem \ref{thm1.3}.

To prove the second part of Theorem \ref{thm1.3}, we consider an initial-boundary value problem for (\ref{1.4})
with initial and boundary data:
\begin{equation}
(\rho,u,b)|_{t=0}=(\overline\rho,0,0),\quad u|_{x=0,1}=0,\quad
b(0,t)=b_1(t)\quad{\rm and}\quad b(1,t)=b_2(t),\label{4.3}
\end{equation}
where $\overline\rho\equiv {\rm const.}>0$ and $b_1(t), b_2(t)$ are
the same as in (\ref{4.1}) and (\ref{1.10}), respectively.  In
the following, for simplicity but without any confusion, we denote
by $(\rho,u,b)$ the solution of the problem (\ref{1.4}), (\ref{4.3}).

\begin{re}\label{re4.1}
 It is worth mentioning that all the
uniform estimates established in Section 3 also hold for the
solution $(\rho,u,b)$ of the problem (\ref{1.4}), (\ref{4.3}).
\end{re}

In order to study the thickness of magnetic boundary layers, we need
to utilize the trivial solution (\ref{4.2}) of the problem
(\ref{1.5}), (\ref{4.1}), the specific initial data in (\ref{4.3})
and the convergence rates in (\ref{1.14}) to obtain the following
improved estimates.
\begin{lem}\label{lem4.1} Let $(\rho,u,b)$ be a smooth solution of (\ref{1.4}), (\ref{4.3}) on
$[0,1]\times[0,T)$. Then there exists a positive constant $C>0$,
independent of $\nu$, such that
\begin{eqnarray}  &&
\sup_{0\leq t< T}\left(\|\rho-\overline\rho\|_{L^2}^2+\|b\|_{L^2}^2+\|u\|_{L^2}^2\right)(t)+\int_0^T\|u_x\|_{L^2}^2dt
\leq C\nu^{1/2},  \label{4.4}    \\
&&  \sup_{0\leq t< T}\|u_x(t)\|_{L^2}^2+\int_0^T\|\rho^{1/2}\dot
u\|_{L^2}^2 dt\leq C\nu^{1/2}.\label{4.5}
\end{eqnarray}
\end{lem}
\pf As an immediate consequence of (\ref{1.14}) and (\ref{4.2}), one obtains
(\ref{4.4}).

To prove (\ref{4.5}), similar to the derivation of (\ref{3.10})$_3$, we get from (\ref{3.7}) and (\ref{3.17}) that
\begin{eqnarray}
&&\frac{\lambda}{2}\frac{d}{dt}\int_0^1u_x^2dx+\int_0^1\rho \dot u^2 dx\nonumber\\
&&\quad=\frac{d}{dt}\int_0^1 \left(P(\rho) + \frac{1}{2} b^2\right)
u_{x} dx+\int_0^1\left(\gamma P(\rho) u_x^2 +b^2u_x^2-\nu b b_{xx}
u_x-\frac{\lambda}{2}
u_x^3  \right)dx\nonumber\\
&&\quad\leq \frac{d}{dt}\int_0^1 \left(P(\rho)-P(\overline\rho) +
\frac{1}{2} b^2\right) u_{x}
dx+C\left(1+\|u_x\|_{L^\infty}\right)\|u_x\|_{L^2}^2 +C\nu^2
\|b_{xx}\|_{L^2}^2,\label{4.6}
\end{eqnarray}
where we have also used the non-slip boundary conditions
$u|_{x=0,1}=0$.

Due to the vanishing initial data given in (\ref{4.3}), we have
\begin{equation}
\left.\int_0^1 \left(P(\rho)-P(\overline\rho) + \frac{1}{2}
b^2\right) u_{x} dx\right|_{t=0}=0,\label{4.7}
\end{equation}
and by (\ref{3.7}), (\ref{3.17}) and (\ref{4.4}), we get
\begin{eqnarray}
\left|\int_0^1 \left(P(\rho)-P(\overline\rho) + \frac{1}{2}
b^2\right) u_{x} dx\right|&\leq&
C\left(\|\rho-\overline\rho\|_{L^2}+\|b\|_{L^2}\right)\|u_x\|_{L^2}\nonumber\\
&\leq&
\frac{\lambda}{4}\|u_x\|_{L^2}^2+C\left(\|\rho-\overline\rho\|_{L^2}^2+\|b\|_{L^2}^2\right)\nonumber\\
&\leq& \frac{\lambda}{4}\|u_x\|_{L^2}^2+C\nu^{1/2}.\label{4.8}
\end{eqnarray}
Thus, by virtue of (\ref{4.7}), (\ref{4.8}) and  (\ref{3.19}) we
deduce after integrating (\ref{4.6}) over $(0,t)$ that for any
$0\leq t< T$,
\begin{equation*}
\|u_x(t)\|_{L^2}^2+\int_0^t\|\rho^{1/2}\dot u\|_{L^2}^2 ds\leq
C\nu^{1/2}+C\int_0^t\|u_x\|_{L^\infty}\|u_x\|_{L^2}^2ds,
\end{equation*}
which, together with (\ref{3.18}), (\ref{3.21}) and Gronwall's inequality, gives (\ref{4.5}).\hfill$\square$

\vskip 2mm

The analysis of the thickness of magnetic boundary layers relies on
the following weighted (interior) $L^2$-estimate of the density and magnetic field.
\begin{lem}\label{lem4.2}
Let $(\rho,u,b)$ be a smooth solution of (\ref{1.4}), (\ref{4.3}) on
$[0,1]\times[0,T)$. Then there exists a positive constant $C>0$,
independent of $\nu$, such that
\begin{equation}
\sup_{0\leq t< T}\int_0^1 \xi(x)\left(|
\rho_x|^2+|b_x|^2\right)(x,t) dx\leq C\nu^{1/2},\label{4.9}
\end{equation}
where $\xi(x)\triangleq x^2(1-x)^2$.
\end{lem}
\pf Differentiating (\ref{1.4})$_3$ with respect to $x$ gives
$$
b_{xt}-\nu b_{xxx}=-u b_{xx}-2 u_xb_x-u_{xx}b,
$$
which, multiplied by $b_x\xi(x)$ in $L^2$ and integrated by parts, yields
\begin{eqnarray}
&&\frac{1}{2}\frac{d}{dt}\int_0^1 b_x^2\xi(x)dx+\nu\int_0^1
b_{xx}^2\xi(x)dx\nonumber\\
&&\quad =\frac{\nu}{2}\int_0^1 b_x^2\xi''(x) dx+\frac{1}{2}\int_0^1u b_x^2\xi'(x)dx\nonumber\\
&&\qquad-\frac{3}{2}\int_0^1 u_x
b_x^2\xi(x) dx-\int_0^1 u_{xx} b b_x\xi(x)dx\triangleq \sum_{i=1}^4I_i,\label{4.10}
\end{eqnarray}
where we have used the fact that $\xi(x)=\xi'(x)=0$ on the boundaries $x=0,1$.

We are now in a position of estimating each term on the right-hand side of (\ref{4.10}). First, it readily follows from (\ref{3.19}) that
\begin{equation}
I_1\leq C\nu\|b_x\|_{L^2}^2\leq C\nu^{1/2}.\label{4.11}
\end{equation}

Noting that due to the non-slip boundary conditions
$u|_{x=0,1}=0$ and (\ref{4.5}) (also cf. (\ref{3.9})), one has
\begin{eqnarray*}
u(x,t)&=& \int_0^x u_x(x,t)dx\leq C x\|u_x\|_{L^\infty},\quad \forall\
x\in[0,1]
\\
u(x,t)&=& -\int_x^1 u_x(x,t)dx\leq C(1-x)\|u_x\|_{L^\infty},\quad\forall\
x\in[0,1],
\end{eqnarray*}
so that, by (\ref{3.18}) and (\ref{3.21}), we have
\begin{eqnarray}
I_2&=&\int_0^1u b_x^2\left[x(1-x)^2-x^2(1-x)\right]dx\nonumber\\
&\leq& C\|u_x\|_{L^\infty}\int_0^1 b_x^2\xi(x)dx\nonumber\\
&\leq& C\left(1+\|\rho^{1/2}\dot u\|_{L^2}\right)\int_0^1 b_x^2\xi(x)dx.\label{4.12}
\end{eqnarray}

Similarly,
\begin{equation}
I_3\leq C\|u_x\|_{L^\infty}\int_0^1 b_x^2\xi(x)dx\leq C\left(1+\|\rho^{1/2}\dot u\|_{L^2}\right)\int_0^1 b_x^2\xi(x)dx.\label{4.13}
\end{equation}

Finally, recalling that
$$
u_{xx}=\lambda^{-1}\left(F_x+P(\rho)_x+bb_x\right),\quad F_x=\rho\dot u ,
$$
and using (\ref{3.7}), (\ref{3.17}) and Cauchy-Schwarz's inequality,
we infer from (\ref{3.8}) that
\begin{eqnarray}
I_4&=&- \lambda^{-1}\int_0^1\left(F_x+P(\rho)_x+bb_x\right)b b_x\xi(x)dx\nonumber\\
&\leq& C\int_0^1\left(|F_x|+|\rho_x|+|b_x|\right)|b_x|\xi(x)dx\nonumber\\
&\leq& C\|F_x\|_{L^2}^2+C\int_0^1\left(\rho_x^2+
b_x^2\right)\xi(x)dx\nonumber\\
&\leq& C\|\rho^{1/2}\dot u\|_{L^2}^2+C\int_0^1\left(\rho_x^2+
b_x^2\right)\xi(x)dx.\label{4.14}
\end{eqnarray}

Thus, inserting (\ref{4.11})--(\ref{4.14}) into (\ref{4.10}), we arrive at
\begin{eqnarray}
&& \frac{d}{dt}\int_0^1 b_x^2\xi(x)dx+\nu\int_0^1
b_{xx}^2\xi(x)dx\nonumber\\
&&\quad\leq C\left(\nu^{1/2}+\|\rho^{1/2}\dot u\|_{L^2}^2\right)+
C\left(1+\|\rho^{1/2}\dot u\|_{L^2}\right)\int_0^1\left(\rho_x^2+
b_x^2\right)\xi(x)dx.\label{4.15}
\end{eqnarray}

In a similar manner, differentiating (\ref{1.4})$_1$ with respect to $x$, multiplying it by $\rho_x\xi(x)$ in $L^2$, and integrating by parts, we deduce
\begin{eqnarray*}
\frac{1}{2}\frac{d}{dt}\int_0^1 \rho_x^2\xi(x)dx
&=&\frac{1}{2}\int_0^1u \rho_x^2\xi'(x)dx-\frac{3}{2}\int_0^1 u_x
\rho_x^2\xi(x) dx-\int_0^1 u_{xx} \rho \rho_x\xi(x)dx\nonumber\\
&\leq&C\|F_x\|_{L^2}^2+C\left(1+\|u_x\|_{L^\infty}\right)\int_0^1\left(\rho_x^2+b_x^2\right)\xi(x)dx\nonumber\\
&\leq& C\|\rho^{1/2}\dot u\|_{L^2}^2+C\left(1+\|\rho^{1/2}\dot u\|_{L^2}\right)\int_0^1\left(\rho_x^2+b_x^2\right)\xi(x)dx,
\end{eqnarray*}
which, combined with (\ref{4.15}), shows that
\begin{eqnarray}
&& \frac{d}{dt}\int_0^1\left(\rho_x^2+ b_x^2\right)\xi(x)dx+\nu\int_0^1
b_{xx}^2\xi(x)dx\nonumber\\
&&\quad\leq C\left(\nu^{1/2}+\|\rho^{1/2}\dot u\|_{L^2}^2\right)+
C\left(1+\|\rho^{1/2}\dot u\|_{L^2}\right)\int_0^1\left(\rho_x^2+
b_x^2\right)\xi(x)dx.\label{4.16}
\end{eqnarray}

Therefore, by virtue of (\ref{4.5}), we deduce from (\ref{4.16}) that
\begin{eqnarray*}
&&\sup_{0\leq t\leq T}\int_0^1\xi(x)\left(\rho_x^2+ b_x^2\right)(x,t)dx+\nu\int_0^T\int_0^1 \xi(x)
b_{xx}^2(x,t)dxdt\nonumber\\
&&\quad \leq  C\left( \nu^{1/2}+\int_0^T\|\rho^{1/2}\dot u\|_{L^2}^2dt\right)\exp\left\{C\int_0^T\left(1+\|\rho^{1/2}\dot u\|_{L^2}\right)dt\right\}\nonumber\\
&&\quad\leq C\nu^{1/2}.
\end{eqnarray*}
The proof of (\ref{4.9}) is complete.\hfill$\square$

\vskip 2mm

With the help of (\ref{4.4}) and (\ref{4.9}), we can estimate the thickness of magnetic boundary layers.

\vskip 2mm

{\it Proof of (ii)-Theorem \ref{thm1.3}.} On one hand, by
(\ref{4.9}) we observe that for any $\delta\in(0,1/2)$,
\begin{eqnarray*}
&&\delta^2\int_\delta^{1-\delta} \left(\rho_x^2+ b_x^2\right)(x,t)dx\nonumber\\
&&\quad=\delta^2\int_\delta^{1/2} \left(\rho_x^2+ b_x^2\right)(x,t)dx
+\delta^2\int_{1/2}^{1-\delta}\left(\rho_x^2+ b_x^2\right)(x,t)dx\nonumber\\
&&\quad\leq \int_\delta^{1/2} x^2\left(\rho_x^2+ b_x^2\right)(x,t)dx+\int_{1/2}^{1-\delta}(1-x)^2\left(\rho_x^2+ b_x^2\right)(x,t)dx\nonumber\\
&&\quad\leq 4\int_\delta^{1/2} x^2(1-x)^2\left(\rho_x^2+ b_x^2\right)(x,t)dx+
4\int_{1/2}^{1-\delta}x^2(1-x)^2\left(\rho_x^2+ b_x^2\right)(x,t)dx\nonumber\\
&&\quad\leq 4\int_\delta^{1-\delta} x^2(1-x)^2\left(\rho_x^2+ b_x^2\right)(x,t)dx\leq C\nu^{1/2},
\end{eqnarray*}
from which it follows that for any $\delta\in(0,1/2)$,
\begin{equation}\label{4.17}
\|\rho_x(t)\|_{L^2(\Omega_\delta)}+\|b_x(t)\|_{L^2(\Omega_\delta)}\leq
C\delta^{-1}\nu^{1/4},\quad\forall\ t\in[0,T),
\end{equation}
where $\Omega_\delta\triangleq(\delta,1-\delta)$. Thus, in view of
(\ref{4.4}) and (\ref{4.17}), we have by Sobolev's inequality that for
any $t\in[0,T)$,
\begin{eqnarray}
&&\|(\rho-\overline\rho)(t)\|_{C(\overline\Omega_\delta)}^2+\|b(t)\|_{C(\overline\Omega_\delta)}^2\nonumber\\
&&\quad\leq C\left(\|(\rho-\overline\rho)(t)\|_{L^2(\Omega_\delta)}^2+\|b(t)\|_{L^2(\Omega_\delta)}^2\right)\nonumber\\
&&\qquad+C\left(\|(\rho-\overline\rho)(t)\|_{L^2(\Omega_\delta)}\|\rho_x(t)\|_{L^2(\Omega_\delta)}
+\|b(t)\|_{L^2(\Omega_\delta)}\|b_x(t)\|_{L^2(\Omega_\delta)}\right)\nonumber\\
&&\quad\leq C\left(\|(\rho-\overline\rho)(t)\|_{L^2(\Omega)}^2+\|b(t)\|_{L^2(\Omega)}^2\right)\nonumber\\
&&\qquad+C\left(\|(\rho-\overline\rho)(t)\|_{L^2(\Omega)}\|\rho_x(t)\|_{L^2(\Omega_\delta)}
+\|b(t)\|_{L^2(\Omega)}\|b_x(t)\|_{L^2(\Omega_\delta)}\right)\nonumber\\
&&\quad\leq C \nu^{1/2}+C\nu^{1/2}\delta^{-1} \leq C\nu^{1/2}\delta^{-1} \to0\quad{\rm as}\quad \nu\to0,\label{4.18}
\end{eqnarray}
provided $\delta=\delta(\nu)$ satisfies
$$
\delta(\nu)\to0\quad{\rm and}\quad \frac{\delta(\nu)}{\nu^{1/2}}\to \infty,\quad{\rm as}\quad\nu\to0.
$$

On the other hand, it is easy to see that
$$
\liminf_{\nu\to0}\left(\|(\rho-\overline\rho)(t)\|_{L^\infty(0,T;C(\overline\Omega))}^2+\|b(t)\|_{L^\infty(0,T;C(\overline\Omega))}^2\right)>0,
$$
provided the boundary data $b_1(t),b_2(t)$ are not identically zero (i.e., $b_1(t),b_2(t)\neq 0$). This, together with (\ref{4.18}), proves the second part of Theorem \ref{thm1.3}.\hfill$\square$

\small{
}
\end{document}